\documentclass[amssymb,amscd,amsfonts,refcheck,11pt,graphicx,verbatim,righttag]{amsart}

\usepackage{amsfonts}
\usepackage{amssymb}

 \numberwithin{equation}{section}
 \allowdisplaybreaks[4]
\theoremstyle{plain}

\newcommand{\N}{{\mathbb N}}

\newcommand{\Z}{{\mathbb Z}}
\newcommand{\R}{{\mathbb R}}

\newcommand{\Q}{{\mathbb Q}}

\newcommand{\G}{{\mathcal G}}
 \newcommand{\NN}{{\mathcal N}}
\newcommand{\D}{{\mathcal D}}

\newcommand{\HH}{{\mathcal H}}
\newcommand{\A}{{\mathcal A}}
\newcommand{\B}{{\mathcal B}}
\newcommand{\M}{{\mathcal M}}

\def\tr{\mathop{\mbox{\rm tr}}\nolimits}

\newtheorem{Thm}{Theorem}
\newtheorem{thm}[Thm]{Theorem}
\newtheorem{cor}{Corollary}

\newtheorem{pro}{Proposition}
\newtheorem{prop}[pro]{Proposition}
\newtheorem{Rem}{Remark}
\newtheorem{rem}[Rem]{Remark}
\newtheorem{Lem}{Lemma}
\newtheorem{lem}[Lem]{Lemma}
\newtheorem{Def}{Definition}
\newtheorem{defi}[Def]{Definition}

\def \proof{\bigbreak\noindent{\it Proof.~~}}

 \numberwithin{equation}{section}

\begin{document}
\title{The spectral properties  of the strongly coupled  Sturm Hamiltonian of eventually constant type}

{\author{Yan-Hui QU}}

\address[Y.-H. QU]{Department of  Mathematics,  Tsinghua University, Beijing 100084, P. R.  China.}
\email{yhqu@math.tsinghua.edu.cn; yanhui.qu@gmail.com}

\begin{abstract}
We study the spectral properties  of the Sturm Hamiltolian of eventually constant type, which includes the Fibonacci Hamiltonian. Let $s$ be the Hausdorff dimension of the spectrum. For $V>20$, we show that the restriction of the $s$-dimensional Hausdorff measure  to the spectrum is a Gibbs type measure; the density of states measure is a Markov measure. Based on the fine structures of these measures,  we show that both measures are exact dimensional; we obtain exact asymptotic behaviors for the optimal H\"older exponent and  the Hausdorff dimension of the density of states measure and for   the Hausdorff dimension of the spectrum. As a consequence, if the frequency is not  silver number type, then for $V$ big enough, we establish strict inequalities between these three spectral characteristics. 

We achieve them by introducing an auxiliary symbolic dynamical system and applying the thermodynamical  and multifractal formalisms of almost additive potentials. 

 \bigskip

\end{abstract}

\maketitle


\section{Introduction}

The Sturm Hamiltonian is a  discrete Schr\"odinger operator
$$
(H_{\alpha,V,\phi}\psi)_n:=\psi_{n-1}+\psi_{n+1}+v_n \psi_n
$$
on $\ell^2(\Z),$ where
the potential $(v_n)_{n\in\mathbb{Z}}$ is given by
\begin{equation}\label{sturm}
v_n=V\chi_{[1-\{\alpha\},1)}(n\alpha+\phi \mod 1),\quad \forall n\in\mathbb{Z},
\end{equation}
where $\alpha>0$ is  irrational, and is called {\it frequency} ($\{\alpha\}$ is the fractional part of $\alpha$),
$V>0$ is called  {\it coupling}, $\phi\in[0,1)$ is
called {\it phase}.  It is well-known that the spectrum of Sturm Hamiltonian
 is independent of   $\phi$,  which  we  denote  by $\Sigma_{\alpha,V}$(see \cite{BIST}).  

Sturm Hamiltonian is firstly  introduced by physicists to model the quasicrystal system, see \cite{BIST} and the references therein  for an excellent  introduction about the physical background. From the mathematical point of view, one is  interested in its spectral properties.

For a discrete Schr\"odinger operator,  several  spectral objects are very important, which are the spectrum, the spectral measure and the density of states measure of the related operator (see \cite{CL} for the exact definitions and related properties).   Let us recall the definition of the density of states measure in our setting. By the spectral theorem, there are Borel probability measures $\mu_{\alpha,V,\phi}$ on $\R$ such that 
$$
\langle \delta_0, g(H_{\alpha,V,\phi})\delta_0\rangle=\int_\R g(x) d\mu_{\alpha,V,\phi}(x)
$$
for all bounded measurable functions $g$, where $\delta_0$ is the element in $\ell^2(\Z)$ which  takes value 1 at site $0$ and  $0$ elsewhere. The {\it density of states measure } $\NN_{\alpha,V}$ is given by the $\phi$-average of these measures with respect to Lebesgue measure, that is,
$$
\int_{\mathbb{T}} \langle \delta_0, g(H_{\alpha,V,\phi})\delta_0\rangle d \phi=\int_\R g(x)d\NN_{\alpha,V}(x)
$$
for all bounded measurable functions $g$. It is well-known that the density of states measure $\NN_{\alpha,V}$ is continuous and supported on $\Sigma_{\alpha,V}$(see for example \cite{CL}).

We will study the fractal  properties  of  $\Sigma_{\alpha,V} $ and $\NN_{\alpha,V}$ for a special class of $\alpha. $  See \cite{F,Fa} for various definitions of fractal dimensions of set and measure.
Write $s_V(\alpha):=\dim_H \Sigma_{\alpha,V}$ and  $d_V(\alpha):=\dim_H \NN_{\alpha,V}$  for the Hausdorff dimensions of $\Sigma_{\alpha,V}$ and $\NN_{\alpha,V}$, respectively; denote by   $\gamma_V(\alpha)$  the optimal H\"older exponent of $\NN_{\alpha,V}$ (see \eqref{optimal-H} for the definition). Before introducing  our results, let us   give a brief survey about the known results on  Sturm Hamiltonian. 
   
 The most prominent model among the Sturm Hamiltonian is the Fibonacci Hamiltonian, for  which the frequency is taken to be  the golden number $\alpha_1=(\sqrt{5}+1)/2.$ This model was  introduced by physicists  to model the quasicrystal  system, see  \cite{KKT,OPRSS}.
S\"ut\"o  \cite{Su} showed that the spectrum  has zero Lebesgue measure for all $V>0$. Then it is  natural  to ask what is the fractal dimension of the spectrum.
It follows implicitly from \cite{Ca} and \cite{Su87} that if $V\ge 16,$ then
$
\dim_B\Sigma_{\alpha_1,V}=\dim_H \Sigma_{\alpha_1,V},
$ where $\dim_B E$ denotes the Box-dimension of $E$ ( see also \cite{DEGT} for explicit statement).
  Raymond \cite{R} first estimated the Hausdorff dimension, he showed that $\dim_H\Sigma_{\alpha_1,V}<1$ for $V>4$.  Jitomirskaya and Last \cite{JL} showed that for any $V>0$, the spectral measure of the operator has positive Hausdorff dimension, as a consequence $\dim_H\Sigma_{\alpha_1,V}>0.$ By using a dynamical method, Damanik et al. \cite{DEGT}   got lower and upper bounds for the dimensions. Due to these bounds they further showed that
\begin{equation}\label{asym-Fibo}
\lim_{V\to\infty} s_{V}(\alpha_1)\ln  V =\ln (1+\sqrt{2}).
\end{equation}
   By studying  the hyperbolicity of  Fibonacci trace map $T_V$ restricting to some invariant surface $S_V$  for small $V$, Damanik and Gorodetski \cite{DG} showed that $s_V(\alpha_1)$ is a $C^\infty$
 function of $V$ on $(0,V_0)$ for some $V_0>0.$ 
 Motivated by totally different problem, Cantat \cite{C} also studied the Fibonacci trace map (indeed more general), he showed that  the $T_V$ is hyperbolic for any $V>0$. Combining with most recent works \cite{P} and \cite{DGY}, now it is known that 
 $s_V(\alpha_1),$ as a function of $V$, is analytic on $(0,\infty)$ and takes values in $(0,1)$.
In \cite{DG2}, Damanik and Gorodetski further showed that $\lim_{V\downarrow 0}s_V(\alpha_1)=1$. In \cite{DG3}, 
Damanik and Gorodetski showed that $\NN_{\alpha_1,V}$ is exact dimensional  and $d_V(\alpha_1)< s_V(\alpha_1)$ for small $V$. In \cite{DG4}, Damanik and Gorodetski showed that $\gamma_V(\alpha_1)\to 1/2$ as $V\to0$ and 
\begin{equation}\label{asym-Fibo-opt}
\lim_{V\to\infty} \gamma_{V}(\alpha_1)\ln  V =\frac{3}{2}\ln \alpha_1.
\end{equation}

Now we  turn to the general Sturm Hamiltonian case. Fix an irrational $\alpha>0$  with continued fraction expansion $[a_0;a_1,a_2,\cdots]$. 
 Bellissard et al.  \cite{BIST} showed
 that $\Sigma_{\alpha,V}$ is a Cantor set of
Lebesgue measure zero.
Damanik, Killip and Lenz \cite{DKL} showed that, if   $\limsup\limits_{k\rightarrow\infty}\frac{1}{k}\sum_{i=1}^k
a_i<\infty$, then  $\dim_H \Sigma_{\alpha,V}>0$. Based on the analysis of Raymond  about  the structure of spectrum \cite{R}, the fractal dimensions of the spectrum of Sturm Hamiltonian were  extensively  studied  in   \cite{LW,LW05,LPW07,FLW,LQW}. 
The current  picture is  the following. Write
$$
K_\ast(\alpha)=\liminf_{k\rightarrow\infty}
(\prod_{i=1}^k a_i)^{1/k}\ \ \ \text{ and }\ \ \  K^\ast(\alpha)=
\limsup_{k\rightarrow\infty}(\prod_{i=1}^k a_i)^{1/k}.
$$
 Fix $V\ge 24$. Then  it is proven in \cite{LW, LQW} that 
$$
\begin{cases}
\dim_H \Sigma_{\alpha,V} \in(0,1) & \text{ if   }\ \  K_\ast(\alpha)< \infty\\
\dim_H \Sigma_{\alpha,V} =1 & \text{ if   }\ \  K_\ast(\alpha)= \infty
\end{cases}
\ \text{ and }\ 
\begin{cases}
\overline{\dim}_B \Sigma_{\alpha,V} \in(0,1) & \text{ if   }\ \  K^\ast(\alpha)< \infty\\
\overline{\dim}_B \Sigma_{\alpha,V} =1 & \text{ if   }\ \  K^\ast(\alpha)= \infty
\end{cases}.
$$
where $\overline{\dim}_B E$ denote the upper Box-dimension of $E$.
Raymond \cite{R}, Liu and Wen \cite{LW} showed that the spectrum $\Sigma_{\alpha,V}$ has a natural covering structure. This structure makes it possible to define the so-called pre-dimensions $s_\ast(\alpha,V)$ and $s^\ast(\alpha,V)$. Then  it is proven in \cite{LPW07, FLW,LQW} that 
$$
\dim_H \Sigma_{\alpha,V}=s_*(\alpha,V)
\ \ \ \text{ and }\ \ \  \overline{\dim}_B \Sigma_{\alpha,V}=s^*(\alpha,V).
$$
 Moreover there exist two constants $0<\rho_\ast(\alpha)\le \rho^\ast(\alpha) $ such that 
\begin{equation}\label{asym-general}
\lim_{V\to \infty} s_*(\alpha,V)\ln V
= \rho_\ast(\alpha)\ \ \ \ \text{ and }\ \ \ \ 
\lim_{V\to \infty} s^*(\alpha,V)\ln V
= \rho^\ast(\alpha).
\end{equation}
 It is proven in  \cite{LQW} that 
$s_*(\alpha,V)$ and $s^*(\alpha,V)$ are
  Lipschitz continuous on any bounded interval of $[24,\infty).$
  
  Recently several works deal with some sub-classes  of Sturm Hamiltonian.  Girand \cite{Gi} considered the frequency $\alpha$ for which the related potential can also be generated  by a primitive invertible substitution.  Mei \cite{Mei} considered the frequency $\alpha$ which has eventually periodic continued fraction expansion,  a strictly larger class than that considered by Girand. In both papers they showed that $\lim_{V\to0}d_V(\alpha)=\lim_{V\to0}s_V(\alpha)=1$, 
   $\NN_{\alpha,V}$ is exact dimensional  and $d_V(\alpha)< s_V(\alpha)$ for small $V$. This generalizes the results in \cite{DG3}.
   Munger \cite{Mu} considered the frequency  of {\it constant type}, i.e. $\alpha_\kappa=[\kappa;\kappa,\kappa,\cdots]$. He gave estimations on the optimal H\"older exponent $\gamma_V(\alpha_\kappa)$ and showed the following asymptotic formula:
   \begin{equation}\label{Munger}
\lim_{V\to\infty} \gamma_V(\alpha_\kappa)\ln V=
\begin{cases}
\frac{3}{2}\ln \alpha_1&\kappa=1\\
\frac{2}{\kappa}\ln \alpha_\kappa&\kappa\ge 2.
\end{cases}
\end{equation}

  In this paper we  will  consider   the frequency of eventually  constant type.  Fix $\kappa\in\N$, define 
  $$
  \mathcal{F}_\kappa:=\{\alpha: \alpha \text{ has expansion } [a_0;a_1,\cdots,a_n,\kappa,\kappa,\cdots]; a_0\ge 0; a_i\in \N, 1\le i\le n; n\in \N\}.
  $$
  Notice that $\alpha_\kappa\in \mathcal{F}_\kappa$ and $\alpha_1\in \mathcal{F}_1$ is the Golden number $(\sqrt{5}+1)/2 $. Define $\mathcal F:=\bigcup_{\kappa=1}^\infty \mathcal{F}_\kappa.$ Any $\alpha\in \mathcal{F}$ is called a frequency of {\it eventually constant type}.
We will study the spectral property of the related Sturm Hamiltonian for large coupling constant $V$.

  Given a probability measure $\mu$ defined on a compact metric space $X$.    Fix $x\in X$, we define the {\it local upper} and {\it lower }dimensions of $\mu$ at $x$ as
$$
\overline{d}_\mu(x)=\limsup_{r\to0}\frac{\ln \mu(B(x,r))}{\ln r}\ \ \ \text{ and }\ \ \ \underline{d}_\mu(x)=\liminf_{r\to0}\frac{\ln \mu(B(x,r))}{\ln r}.
$$
In the case  $\overline{d}_\mu(x)=\underline{d}_\mu(x)$, we say that the {\it local  dimension} of   $\mu$ at $x$ exists and we denote it by  $d_\mu(x)$. 
The {\it optimal H\"older exponent } $\gamma_\mu$ of $\mu$ is defined  as 
\begin{equation}\label{optimal-H}
\gamma_\mu:=\inf\{\underline{d}_\mu(x):x\in X\}.
\end{equation}
The Hausdorff dimension of $\mu$ is defined as
  \begin{equation}\label{dim-meas}
  \dim_H\mu:=\sup\{s: \underline{d}_\mu(x)\ge s \text{ for  } \mu \text{ a.e. }x\in X\}.
  \end{equation}

If there exists a constant $d$ such that $d_\mu(x)=d$ for $\mu$ a.e. $x\in X$, then necessarily $\dim_H\mu=d$. In this case we call $\mu$ {\it exact dimensional}.

 Our main result is the following.

\begin{thm}\label{main} 
Fix $V>20$ and $\kappa\in \N.$ Then

(i)\ There exist three positive numbers $\gamma_V(\kappa), d_V(\kappa)$ and $s_V(\kappa)$ such that  for any $\alpha\in \mathcal{F}_\kappa$, 
\begin{equation}\label{uniform-dim}
\gamma_V(\alpha)=\gamma_V(\kappa), \ \ \ d_V(\alpha)= d_V(\kappa) \ \ \ \text{ and } \ \ \ s_V(\alpha)=s_V(\kappa).
\end{equation}

Moreover for any fixed $\alpha\in \mathcal{F}_\kappa$,  the following two assertions hold: 

(ii)\ $\HH^{s_V(\kappa)}|_{\Sigma_{\alpha,V}}$ is a Gibbs type measure. Consequently, $\HH^{s_V(\kappa)}|_{\Sigma_{\alpha,V}}$ is exact dimensional and  $0<\HH^{s_V(\kappa)}(\Sigma_{\alpha,V})<\infty.$

(iii)\ $\NN_{\alpha,V}$ is a Markov  measure and  a Gibbs type measure. Consequently, it is exact dimensional. 

(iv)\ For each $\kappa\in\N$, there exist three constants $0<\hat \varrho_\kappa\le \varrho_\kappa\le \rho_\kappa$ such that 
$$
\lim_{V\to\infty} \gamma_V(\kappa) \ln V=\hat \varrho_\kappa,\ \ \lim_{V\to\infty} d_V(\kappa) \ln V=\varrho_\kappa\ \ \text{ and }\ \ \lim_{V\to\infty} s_V(\kappa) \ln V=\rho_\kappa.
$$
Moreover,    $\hat \varrho_2=\varrho_2=\rho_2=\ln (1+\sqrt{2})$ and $\hat \varrho_\kappa<\varrho_\kappa<\rho_\kappa$ when $\kappa\ne 2.$

(v)\ 
When  $\kappa\ne 2$, there exists $V_0(\kappa)>20$ such that  for all $V\ge V_0(\kappa)$ we have 
$$
\gamma_V(\kappa)<d_V(\kappa)<s_V(\kappa).
$$
 \end{thm}
 
 \begin{rem} \label{rem-main}
 {\rm
 
 (1)  \eqref{uniform-dim} shows that the three quantities $\gamma_V(\alpha), d_V(\alpha)$ and $s_V(\alpha)$ only depend on the ``tail" of the expansion of $\alpha$ for $\alpha$ of eventually constant type. Although this is expected intuitively, the proof is far from trivial. 
 
  (2)  See Theorem \ref{gibbs} (i) for the definition of Gibbs measure.
 See Definition \ref{gibbs-type}(Section \ref{HDS}) for the formal definition of Gibbs type measure. Roughly speaking 
  ``Gibbs type measure" means that  it  can be decomposed to finite pieces such that each piece is strongly equivalent to an image of a Gibbs measure under a bi-Lipschitz map (see Section \ref{thermo-multi} for the definition of strongly equivalence of measures). The assertions  (ii) and (iii) tell us  that both  measures $\HH^{s_V(\kappa)}|_{\Sigma_{\alpha,V}}$ and $\NN_{\alpha,V}$ have very good dynamical structure. 
  
  (3)    For $\alpha=\alpha_\kappa$, the first equality in (iv) is   \eqref{Munger} ( see \cite{Mu}); for any irrational $\alpha,$ the third equality in (iv) is \eqref{asym-general} (see  \cite{LPW07,FLW,LQW}). We state them  here for  comparison. Our method gives a new proof for \eqref{Munger}.
 
 (4) Frequency $\alpha_1$ corresponds to the   Fibonacci Hamiltonian.  By Remark \ref{hat-rem-varrho}, \ref{rem-varrho} and \ref{rem-rho}  in Section \ref{APC}, we have 
 \begin{equation}\label{two-relation}
\hat \varrho_1=\frac{3}{2}\ln\alpha_1,\ \ \varrho_1=\frac{5+\sqrt{5}}{4}\ln \alpha_1\ \ \text{ and }\ \ \rho_1=\ln (1+\sqrt{2}).
\end{equation}
The first and the third equalities   in  \eqref{two-relation} are known, which are \eqref{asym-Fibo-opt} and \eqref{asym-Fibo},  respectively, see \cite{DEGT,DG4}. But the second one in \eqref{two-relation}  is new. 

(5) Recall that in \cite{DG3}, $d_V(1)<s_V(1)$ is proven for $V$ small. 
Here for any $\kappa\ne 2$ and $V$ is large, we show that $d_V(\kappa)<s_V(\kappa)$. As explained in \cite{DG3}, $\NN_{\alpha,V}$ is the harmonic measure determined by $\Sigma_{\alpha,V}$. It is a general belief that if a planar set $E$ is dynamically defined, then the Hausdorff  dimension of the harmonic measure determined by $E$ is  strictly less than the Hausdorff dimension of $E$(see for example the survey paper \cite{Ma} or the book \cite{GM}).
Our results verify this belief  in this special situation. 

(6) In the case   $\kappa=2$,  $\alpha_2$ is sometimes called {\it silver } number. We call $\alpha\in \mathcal{F}_2$ a {\it silver type} number.  By  the explanation above, we still expect $d_V(2)<s_V(2)$. However,  since $\varrho_2=\rho_2$, we can not achieve this by the asymptotic formulas anymore. Finer estimations are needed. This makes the silver type  number case an interesting object to study. We note that in \cite{DM}, the trace map related to silver number has been studied. They  showed that the non-wandering set of this map is hyperbolic if the coupling is sufficiently large.  We also remark that, it follows from the general theories developed in \cite{C} that  the non-wandering set  is hyperbolic for all coupling constants.

(7) Let us  say a few words about the proof. Based  on the analysis in  \cite{R,LW, FLW} about the nested structure of the spectrum, we can introduce an auxiliary  symbolic dynamical system which codes certain subsets of   the spectrum.  By introducing  two potentials (one is additive, which is related to the density of states measure $\NN_{\alpha,V}$; the other  is almost additive, which is related to the Hausdorff measure restricted to the spectrum ${\Sigma_{\alpha,V}}$), we can use the thermodynamical and multifractal formalisms to analyze the spectrum and $\NN_{\alpha,V}$,  respectively. We  show that both measures have the  Gibbs property, moreover $\HH^{s_V(\alpha)}|_{\Sigma_{\alpha,V}}$ is kind of measure of maximal dimension (see Section \ref{HDS}, Remark \ref{maximal-dim} and Theorem \ref{main-Hausdorff}) and $\NN_{\alpha,V}$ is kind of measure of maximal entropy (see Section \ref{DSMO}, Theorem \ref{max-entropy} and Theorem  \ref{main-DSM}). These good structures in turn give very exact and explicit informations  about the dimensions of the spectrum and the density of states measure. 
}
 \end{rem}

Shortly after we finished the first version of  our paper, we saw the   impressive paper \cite{DGY}.  In \cite{DGY}, Damanik, Gorodetski and Yessen  completed the picture for the Fibonacci Hamiltonian (the case $\alpha=\alpha_1$) by getting rid of the smallness and largeness restriction on the coupling constant $V$, giving the explicit formulas for $\gamma_V(\alpha_1), d_V(\alpha_1), s_V(\alpha_1)$ (and another important quantity-- transport exponent) and the exact asymptotic behaviors of these quantities(among many other things). It seems interesting to  make some comments on their methods and point out some connections of their work with ours.  In their paper, they consider the so-called Fibonacci trace map $T_V$ and the related maximal hyperbolic invariant set $\Lambda_V$.  It is well known that the spectrum is obtained by intersecting  a special line $\ell_V$ with the stable lamination of  $\Lambda_V$.   Through the previous works of the authors \cite{DG,DG2, DG3}, to obtain the spectral  properties of the spectrum, a crucial step is  to show that $\ell_V$ intersects transversally with the stable lamination  of $\Lambda_V$ for any $V>0.$ In  \cite{DGY}, they  made  a decisive  progress and   succeeded  to achieve this step. Then by using powerful tools from hyperbolic dynamical system, they derived the complete picture for Fibonacci Hamiltonian. Now let us point out some connections of their work with ours.    In \cite{DGY} Theorem 1.6, they got  dimension formulas for $\Sigma_{\alpha_1,V}$ (formula (10)) and $\NN_{\alpha_1,V}$ (formula (11)) and formula for the optimal H\"older exponent of $\NN_{\alpha_1,V}$ (formula (12)).  In our setting the counterparts are \eqref{dim-maxi-haus}, \eqref{NN-dim-for} and  \eqref{formula-optimal}. In \cite{DGY} Theorem 1.10, they got exact asymptotic formulas for $\gamma_V(\alpha_1), d_V(\alpha_1), s_V(\alpha_1)$. In our case this  is related  to Remark \ref{rem-main} (4). 
  
  The rest of the paper are organized as follows. In Section \ref{prepare},
 we introduce the notations and summarize known results and connections which we will use. In Section \ref{geo-lemma}, we prove a geometric lemma, which will be used to establish a bi-Lipschitz equivalence between the dynamical subset of  the  spectrum and the symbolic space. In Section \ref{HDS}, we study the Hausdorff dimension and Hausdorff measure of a dynamical subset of the spectrum. In Section \ref{DSMO}, we study the dimension properties of the density of states measure $\NN_{\alpha,V}$ restricted to a dynamical subset. In Section \ref{MAOE}, we conduct the multifractal analysis of $\NN_{\alpha,V}$, in particular,  we get an expression for the optimal H\"older exponent of $\NN_{\alpha,V}$. In Section \ref{global}, by comparing  two arbitrary  dynamical subsets, we obtain the global picture. In particular,  we prove Theorem \ref{main} (i), (ii) and (iii). In Section \ref{APC}, we study the asymptotic properties  and prove Theorem \ref{main} (iv) and (v). Finally in Section \ref{appendix}, we give another proof of the fact that $d_V(\kappa)<s_V(\kappa)$ when $\kappa\ne 2$. 

\section{Preliminaries }\label{prepare}

In this section we summarize  known results and connections which we will use.

At first we discuss the structure of the spectrum and  give a coding for it. Next we collect some useful facts about  Sturm Hamiltonian. Finally  
we recall the thermodynamical and multifractal formalisms for the almost additive potentials.

\subsection{The structure of the spectrum }\

We describe  the structure of the spectrum
$\Sigma_{\alpha,V}$ for fixed  $\alpha\in \R\setminus\Q$ and $V>0$. 
We also collect some known facts that will be used later, for more
details, we refer to \cite{BIST,R,LW,T}.

Since $\Sigma_{\alpha,V}$ is independent with the phase $\phi$, in the rest of the paper we can and will take $\phi=0.$ Assume $\alpha$ has continued fraction expansion $[a_0;a_1,a_2,\cdots]. $
Let   $p_n/q_n (n>0)$ be the $n$-th partial quotient  of $\{\alpha\}$, the fractional part of $\alpha$, 
given by:
\begin{eqnarray}
\nonumber &&p_{-1}=1,\quad p_0=0,\quad p_{n+1}=a_{n+1} p_n+p_{n-1},\ n\ge 0,\\
\label{recur-q} &&q_{-1}=0,\quad q_0=1,\quad q_{n+1}=a_{n+1} q_n+q_{n-1},\ n\ge 0.
\end{eqnarray}
 
Let $n\geq1$ and $x\in\mathbb{R}$, the transfer matrix $M_n(x)$
over $q_n$ sites is defined by
$${\mathbf M}_n(x):=
\left[\begin{array}{cc}x-v_{q_n}&-1\\ 1&0\end{array}\right]
\left[\begin{array}{cc}x-v_{q_n-1}&-1\\ 1&0\end{array}\right]
\cdots
\left[\begin{array}{cc}x-v_1&-1\\ 1&0\end{array}\right],
$$
where $v_n$ is defined in \eqref{sturm}. By convention we  take
$$
\begin{array}{l}
{\mathbf M}_{-1}(x)= \left[\begin{array}{cc}1&-V\\
0&1\end{array}\right]
\end{array}
\ \ \ \text{ and }\ \ \ 
\begin{array}{l}
{\mathbf M}_{0}(x)= \left[\begin{array}{cc}x&-1\\
1&0\end{array}\right].
\end{array}
$$

\smallskip

For $n\ge0$, $p\ge-1$, let $t_{(n,p)}(x)=\tr {\mathbf M}_{n-1}(x) {\mathbf M}_n^p(x)$ and
$$
\sigma_{(n,p)}=\{x\in\mathbb{R}:|t_{(n,p)}(x)|\leq2\},
$$
 where  $\tr M$
stands for the trace of the matrix $M$.
Then
$$
(\sigma_{(n+2,0)}\cup\sigma_{(n+1,0)})\subset
(\sigma_{(n+1,0)}\cup\sigma_{(n,0)}).
$$
Moreover
$$\Sigma_{\alpha,V}=\bigcap_{n\ge0}(\sigma_{(n+1,0)}\cup\sigma_{(n,0)}).$$

The intervals in  $\sigma_{(n,p)}$ are called  {\em bands}. Take some band 
$B\in \sigma_{(n,p)}$,  then $t_{(n,p)}(x)$ is monotone on
$B$ and
$t_{(n,p)}(B)=[-2,2].$
We call $t_{(n,p)}$ the {\em generating polynomial} of $B$ and denote it by $h_B:=t_{(n,p)}$.

$\{\sigma_{(n+1,0)}\cup\sigma_{(n,0)}:n\ge 0\}$ form a covering of $\Sigma_{\alpha,V}$.
However there are some repetitions between $\sigma_{(n,0)}\cup\sigma_{(n-1,0)}$
and $\sigma_{(n+1,0)}\cup\sigma_{(n,0)}$.
It is possible to choose a covering of $\Sigma_{\alpha,V}$ elaborately such that
we can get rid of these repetitions, as we will describe in the follows:

\begin{defi}{\rm (\cite{R,LW})}
For $V>4$, $n\ge0$, we define three types of bands as follows:

$(n,{\rm I})$-type band: a band of $\sigma_{(n,1)}$ contained in a
band of $\sigma_{(n,0)}$;

$(n,{\rm II})$-type band: a band of $\sigma_{(n+1,0)}$ contained
in a band of $\sigma_{(n,-1)}$;

$(n,{\rm III})$-type band: a band of $\sigma_{(n+1,0)}$ contained
in a band of $\sigma_{(n,0)}$.
\end{defi}

All three types of bands actually occur and they are disjoint. We call these bands {\em
spectral generating bands of order $n$}. Note that  there 
are only two spectral generating bands of order $0$, 
one is $\sigma_{(0,1)}=[V-2,V+2]$ with  generating polynomial $t_{(0,1)}=x-V$ and type $(0,{\rm I})$, the other is  $\sigma_{(1,0)}=[-2,2]$ with generating polynomial  $t_{(1,0)}=x$ and type $(0,{\rm III})$. They are
contained in $\sigma_{(0,0)}=(-\infty,+\infty)$ with
generating polynomial $t_{(0,0)}\equiv2$. For  convenience, we
call $\sigma_{(0,0)}$ the spectral generating band of order $-1$.

For any $n\ge-1$, denote by $\B_n$ the set of spectral
generating bands of order $n$, then the intervals in $\B_n$ are disjoint.
Moreover (\cite{R,LW})
\begin{itemize}
\item $(\sigma_{(n+2,0)}\cup\sigma_{(n+1,0)})\subset
\bigcup_{B\in\B_n}B
\subset (\sigma_{(n+1,0)}\cup\sigma_{(n,0)})$,
thus
$$
\Sigma_{\alpha,V}=\bigcap_{n\ge0}
\bigcup_{B\in\B_n}B;
$$
\item
any $(n,I)$-type band contains only one band in $\B_{n+1}$, which is of  $(n+1,II)$-type.
\item
any $(n,II)$-type band contains $2a_{n+1}+1$ bands in $\B_{n+1}$,
$a_{n+1}+1$ of which are of  $(n+1,I)$-type and $a_{n+1}$ of which are of  $(n+1,III)$-type.
\item
any $(n,III)$-type band contains $2a_{n+1}-1$ bands in $\B_{n+1}$,
$a_{n+1}$ of which are of  $(n+1,I)$-type and $a_{n+1}-1$ of which are of  $(n+1,III)$-type.
\end{itemize}

Thus $\{\B_n\}_{n\ge0}$ forms a natural
covering(\cite{LW05,LPW07}) of the spectrum $\Sigma_{\alpha,V}$.

\subsection{The coding of the spectrum}\label{coding}\

In the following we give a coding of  the spectrum $\Sigma_{\alpha,V}$ based on \cite{LQW},  however, we note that   the presentation  here  is slightly different from that in  \cite{LQW}.

 For each $N\in\N,$ we define an alphabet $\A_N$ as 
$$
\A_N:=\{(I,j)_N:j=1,\cdots, N+1\}\cup \{(II,1)_N\}\cup \{(III,j)_N:j=1,\cdots, N\}.
$$
Then $\#\A_N=2N+2. $
We order the elements in $\A_N$ as 
$$
(I,1)_N<\cdots <(I,N+1)_N<(II,1)_N<(III,1)_N<\cdots< (III,N)_N.
$$
To simplify the notation we rename  the above line as 
$$
e_{N,1}<e_{N,2}<\cdots< e_{N,2N+2}.
$$
Given $e_{N,i}\in\A_N$ and $e_{M,j}\in \A_M$, we call $e_{N,i}e_{M,j}$  {\it admissible}, denote by $e_{N,i}\to e_{M,j},$ if 
\begin{eqnarray*}
(e_{N,i},e_{M,j})&\in& \{((I,k)_N,(II,1)_M): 1\le k\le N+1\}\cup \\
&&\{((II,1)_N,(I,l)_M): 1\le l\le M+1\}\cup \\
&&\{((II,1)_N,(III,l)_M): 1\le l\le M\}\cup\\
&&\{((III,k)_N,(I,l)_M):1\le k\le N,\ 1\le l\le M\}\cup \\
&&\{((III,k)_N,(III,l)_M):1\le k\le N, 1\le l\le M-1\}.
\end{eqnarray*}
For pair $(\A_N, \A_M),$ we define the incidence matrix $A_{NM}= (a_{ij})$ of size $(2N+2)\times(2M+2)$ as 
$$
a_{ij}=
\begin{cases}
1& e_{N,i}\to e_{M,j}\\
0& \text{ otherwise}
\end{cases}
$$
When $N=M,$ we write $A_N:=A_{NN}.$
For any $N\in \N$, 
we define a related matrix as follows:
$$
\hat A_N=
\left[\begin{array}{ccc}
0&1&0\\
N+1&0&N\\
N&0&N-1
\end{array}\right]
$$
 
We also define $\A_0:=\{I, III\}$.  Write $e_1=I$ and $e_2=III$.  Given $e_i\in \A_0$ and $e_{M,j}\in \A_{M}$, we call $e_ie_{M,j}$  {\it admissible}, denote by $e_i\to e_{M,j},$ if 
\begin{eqnarray*}
(e_i,e_{M,j})&\in& \{(I,(II,1)_M)\}\cup 
\{(III,(I,l)_M):1\le l\le M\}\cup \\
&&\{(III,(III,l)_M):1\le l\le M-1\}.
\end{eqnarray*}
For pair $(\A_0, \A_M),$ we define the incidence matrix $A_{0M}= (a_{ij})$ of size $2\times(2M+2)$ as 
$$
a_{ij}=
\begin{cases}
1& e_{i}\to e_{M,j}\\
0& \text{ otherwise}
\end{cases}
$$

Recall that $\alpha=[a_0;a_1,a_2,\cdots]$. Define a symbolic space  $\Omega^{(\alpha)}$   with alphabet sequence  $\{\A_0\}\cup\{ \A_{a_i}: i\ge 1\}$ and incidence matrix sequence  $\{A_{0a_1}\}\cup\{A_{a_ia_{i+1}}: i\ge 1\}$ as 
$$
\Omega^{(\alpha)}=\{e_{i_0}e_{a_1,i_1}e_{a_2,i_2}\cdots: \ \ e_{i_0}\in \A_0; e_{a_j, i_j}\in \A_{a_j}; e_{i_0}\to e_{a_1,i_1}, e_{a_j,i_j}\to e_{a_{j+1},i_{j+1}}, j\ge 1\}.
$$
For $\omega\in \Omega^{(\alpha)}$, we write $\omega|_n=\omega_0\cdots \omega_n.$ More generally we write $\omega[n,\cdots, n+k]$ for $\omega_n\cdots \omega_{n+k}.$ 
Define $\Omega^{(\alpha)}_{n}:=\{\omega|_n: \omega\in \Omega^{(\alpha)}\}$ and $\Omega^{(\alpha)}_{\ast}=\bigcup_{n} \Omega^{(\alpha)}_{n}.$   Given $w=w_0\cdots w_n\in \Omega^{(\alpha)}_{n}$, denote by $|w|$ the length of $w$, then $|w|=n+1$.  If $w_n=(T,j)_{a_n}$, then we write $t_w=T,$ and call $w$ has {\it  type } $T.$  If $w=uw^\prime$, then we say $u$ is a {\it prefix } of $w$ and denote by $u\prec w.$ Given $u=u_0\cdots u_n\in  \Omega^{(\alpha)}_{n}$ and $v=v_0\cdots v_m\in \prod_{j=n}^{n+m}\A_{a_j}$, if $u_n=v_0$, then we write 
 $u\star v:=u_0\cdots u_n v_1\cdots v_m.$ For $\omega,\tilde \omega\in \Omega^{(\alpha)}$, we denote by $\omega\wedge \tilde \omega$ the maximal common prefix of $\omega$ and $\tilde\omega.$
Given $w\in \Omega^{(\alpha)}_{n},$
we define the cylinder 
$$
[w]:=\{\omega\in \Omega^{(\alpha)}: \omega|_n=w\}.
$$

 In the following we explain that $\Omega^{(\alpha)}$ is a coding of the spectrum $\Sigma_{\alpha,V}$. 
 Define 
$B_{I}$ to be the unique $(0,I)$ type band in $\B_0$ and  $B_{III}$ to be the unique $(0,III)$ type band in $\B_0$. Then 
$$
\B_0=\{B_w: w\in \Omega_0^{(\alpha)}\}.
$$
Assume $B_w$ is defined for any $w\in  \Omega_{n-1}^{(\alpha)}$.
Now for any $w\in\Omega_n^{(\alpha)}$, write $w=w^\prime e=w^\prime(T,j)_{a_n}$. Then 
define $B_w$ to be the unique $j$-th  $(n,T)$ type band in $\B_n$ which is contained in $B_{w^\prime}$.  Then 
$$
\B_n=\{B_w: w\in\Omega_n^{(\alpha)}\}.
$$
 Thus we can define a natural projection $\pi:\Omega^{(\alpha)}\to \Sigma_{\alpha,V}$ as 
$$
\pi(\omega):=\bigcap_{n\ge 0} B_{\omega|_n}.
$$ 
 It is seen that $\pi$ is a one-to-one map. So $\Omega^{(\alpha)}$  is a coding of $\Sigma_{\alpha,V}$.
 Write $X_w=\pi([w])$ for 
each $w\in\Omega_\ast^{(\alpha)}$ and denote by $h_w$ the generating polynomial of $B_w.$

For any $n\ge 1$, we also define the following symbolic space:
$$
\Omega^{(\alpha, n)}=\{e_{a_n,i_n}e_{a_{n+1},i_{n+1}}\cdots: \ \  e_{a_j, i_j}\in \A_{a_j};   e_{a_j,i_j}\to e_{a_{j+1},i_{j+1}}, j\ge n\}.
$$
Formally, we obtain  $\Omega^{(\alpha,n)}$ by shifting $\Omega^{(\alpha)}$ to the left $n$ times. In general, all the $\Omega^{(\alpha,n)}$ are different, thus it is hard to define a dynamic on them. However, if $\alpha$ is of eventually constant type and has continued fraction expansion $[a_0;a_1,\cdots, a_{\hat n},\kappa,\kappa,\cdots]$, then it is easy to see that 
all the $\{\Omega^{(\alpha,n)}: n> \hat n\}$ are the same and indeed they are a subshift of finite type, with alphabet $\A_\kappa$ and incidence matrix $A_\kappa=A_{\kappa\kappa}$. For this reason,  we give the following definition:

\begin{defi}
For any $\kappa\in \N$, we denote the subshift of finite type with alphabet $\A_\kappa$ and incidence matrix $A_\kappa$ by  $\Omega^{(\kappa)}$. Together with the shift map $\sigma$, $(\Omega^{(\kappa)}, \sigma)$ becomes a dynamical system. 
\end{defi}

By endowed with the usual metric on $\Omega^{(\kappa)}$, $(\Omega^{(\kappa)}, \sigma)$ becomes a topological dynamical system.  The following lemma implies that indeed $(\Omega^{(\kappa)}, \sigma)$ is topologically mixing.

Recall that a nonnegative square matrix $B$ is called {\it primitive} if there exists some $k\in \N$ such that all the entries of  $B^k$ are positive.

\begin{lem} \label{primitive}
For each $\kappa\in\N$,  $A_\kappa, \hat A_\kappa$ are  primitive and have the same Perron-Frobenius eigenvalue $\alpha_\kappa.$  Moreover 
$\hat A_\kappa$ has three different eigenvalues $\alpha_\kappa, -1,-1/\alpha_\kappa$, consequently $\hat A_\kappa$ is  diagonalizable. 
\end{lem}

\proof\
It is straightforward  to check that $\hat A_\kappa^5>0,$ thus $\hat A_\kappa$ is primitive. 
 We also have 
 $$
 \det(\lambda I_3-\hat A_\kappa)=(\lambda^2-\kappa\lambda-1)(\lambda+1).
 $$
 Thus the Perron-Frobenius eigenvalue  of $\hat A_\kappa$ is $\alpha_\kappa$ and   $\hat A_\kappa$ has three different eigenvalues $\alpha_\kappa, -1,-1/\alpha_\kappa$.
 
 On the other hand if we consider the graph related to the incidence matrix $A_\kappa$, then it is easy  to show that the graph is aperiodic. Consequently $A_\kappa$ is primitive. By direct computation we get that 
 $$
 \det(\lambda I_{2\kappa+2}- A_\kappa)=\lambda^{2\kappa-1}(\lambda^2-\kappa\lambda-1)(\lambda+1).
 $$
 Thus the Perron-Frobenius eigenvalue of $A_\kappa$ is also $\alpha_\kappa.$
\hfill $\Box$

In what follows, we will make essential use of this symbolic dynamic $(\Omega^{(\kappa)},\sigma)$ to understand the structure of the spectrum and the density of states measure of  Sturm Hamiltonian with frequencies of eventually constant type.
 
\subsection{Useful results for Sturm Hamiltonian }\

In this subsection we collect some useful results for Sturm Hamiltonian.

Fix an irrational frequency $\alpha$ with continued fraction expansion $[a_0;a_1,a_2,\cdots]$ and consider the operator $H_{\alpha,V,0}$.
We write $h_k(x):=t_{k+1,0}(x)={\tr M_k(x)}$ and write
\begin{equation}\label{sigmak}
\sigma_k:=\sigma_{(k+1,0)}=\{x\in\R : h_k(x)|\le 2\}.
\end{equation}

The following two lemmas are essentially  proven in \cite{R}:

\begin{lem}\label{sigma-k}
$\sigma_k=\{B\in \B_k: B \text{ is of  type } (k,II) \text{ or } (k,III)\}$. 
\end{lem}

Define $v^I=(1,0,0), v^{II}=(0,1,0), v^{III}=(0,0,1)$ and $v_\ast=(0,1,1)^T.$

\begin{lem}
Given $w\in  \Omega_m^{(\alpha)}$,  then 
\begin{equation}\label{num}
\#\{u\in \Omega_{m+n}^{(\alpha)}:w\prec u, t_u=II,III\}=v^{t_w}\cdot \hat A_{a_{m+1}}\cdots \hat A_{a_{m+n}} \cdot v_\ast.
\end{equation}
\end{lem}

The following theorem is \cite{FLW} Theorem 2.1, see also \cite{LQW} Theorem 3.1.
\begin{thm}[Bounded variation]\label{bvar}
Let $V>20$ and $\alpha=[a_0;a_1,a_2,\cdots]$ be irrational with $a_n$ bounded by $M.$ Then there exists a constant
$\xi=\xi(V,M)>1$  such that for any spectral
generating band $B$   with generating
polynomial $h$,
$$\xi^{-1}\le \left|\frac{h'(x_1)}{h'(x_2)}\right|\le \xi,\quad  \forall x_1,x_2\in B.$$

\end{thm}

The following theorem is \cite{FLW} Theorem 5.1, see also \cite{LQW} Theorem 3.3.
\begin{thm}[Bounded covariation]\label{bco}
Let $V>20$, $\alpha=[a_0;a_1,a_2,\cdots], \tilde \alpha=[\tilde a_0;\tilde a_1,\tilde a_2,\cdots]$ be irrational with $a_n, \tilde a_n$ bounded by $M$. Then there exists 
constant $\eta=\eta(V,M)>1$   such
that if $w, wu\in  \Omega_\ast^{(\alpha)}$ and $\tilde w, \tilde wu\in  \Omega_\ast^{(\tilde \alpha)}$,
then 
$$\eta^{-1} \frac{|{B}_{wu}|}{|{B}_{w}|}\le
\frac{|B_{\widetilde wu}|}{|B_{\widetilde w}|}\le \eta
\frac{|{B}_{wu}|}{|{B}_{w}|}.$$
\end{thm}

We remark that  in \cite{FLW}, they only considered the case $\alpha=\tilde \alpha.$ However by checking the proof, the same argument indeed shows the  stronger result as stated in Theorem \ref{bco}.  

The following lemma is a direct consequence of \cite{FLW} Corollary 3.1:
\begin{lem}\label{band-ratio}
Let $V>20$ and $\alpha=[a_0;a_1,a_2,\cdots]$ be irrational with $a_n$ bounded by $M$. 
Then there exist  constants $0<c_1=c_1(V,M)<c_2=c_2(V,M)<1$ and $n_0=n_0(V,M)\in \N$ such that for any $n>n_0$ and any  $w\in  \Omega_n^{(\alpha)}$ 
$$
c_1|B_{w|_{n-n_0}}|\le|B_{w}|\le c_2|B_{w|_{n-n_0}}|.
$$
There exists  constant $0<c_3=c_3(V,M)<1$ such that for any $w\in  \Omega_n^{(\alpha)}$
\begin{equation}\label{band-length}
c_3^n\le |B_w|\le 2^{2-n}.
\end{equation}
\end{lem}
 
The following lemma is \cite{LQW} Lemma 3.7 by taking $a_i=\kappa.$

 \begin{lem}\label{lem-bc}
 Assume $V>20$ and $\alpha=\alpha_\kappa.$
 Write $t_1=(V-8)/3$ and $t_2=2(V+5).$ For $e_{\kappa,\kappa+2}=(II,1)_\kappa$ and $w=w_0\cdots w_n \in\Omega_n^{(\alpha_\kappa)},$   write $|w|_{e_{\kappa,\kappa+2}}:=\#\{1\le i\le n: w_i= e_{\kappa,\kappa+2}\}$. 
Then 
$$
t_2^{(1-\kappa) |w|_{e_{\kappa,\kappa+2}}}\left(t_2\kappa^3\right)^{|w|_{e_{\kappa,\kappa+2}}-n}\le |B_w|\le 4 t_1^{(1-\kappa) |w|_{e_{\kappa,\kappa+2}}}\left(t_1\kappa\right)^{|w|_{e_{\kappa,\kappa+2}}-n}.
$$
\end{lem}

\subsection{Recall on thermodynamical formalism and multifractal analysis}\label{thermo-multi}\

If $X$ is a compact metric space and $T:X\to X$ is continuous, then $(X,T)$ is called a {\it topological dynamical system}, TDS for short. $\M(X)$ is the set of all probability measures supported on $X.$ $\M(X,T)$ is the set of all $T$-invariant probability measures supported on $X.$ Given $\mu,\nu\in \M(X)$, if there exists a constant $C>1$ such that $C^{-1}\nu\le \mu\le C\nu$, then we say $\mu$ and $\nu$ are {\it strongly  equivalent} and denote by $\mu\asymp \nu.$

 Assume $(X,T)$ is a TDS and  $\Phi=\{\phi_n:n\ge 0\}$ is a family of continuous functions from $X$ to $\R$. If there exists a constant $C(\Phi)\ge0$ such that for any $n,k\ge 0$
$$
| \phi_{n+k}(x)- \phi_n(x)-\phi_k(T^nx)|\le C(\Phi),
$$
then $\Phi$ is called a family of {\it almost additive} potentials. We use $C_{aa}(X,T)$ to denote the set of all the almost additive potentials defined on $X.$ When $C(\Phi)=0$,  it is seen that $\phi_n(x)=\sum_{j=0}^{n-1}\phi_1(T^jx)=:S_n\phi_1(x).$ That is, $\phi_n$ is the {\it ergodic sum } of $\phi_1$. In this case we say $\Phi$ is a family of {\it additive} potentials. 

Given $\Phi\in C_{aa}(X,T)$.  If there exists a constant $c>0$ such that $\phi_n(x)\ge cn$ for any $n\ge 0$, then we say $\Phi$ is {\it positive} and write $\Phi\in C_{aa}^+(X,T).$ Similarly, if there exists a constant $c<0$ such that $\phi_n(x)\le cn$ for any $n\ge 0$, then we say $\Phi$ is {\it negative} and write $\Phi\in C_{aa}^-(X,T).$

 Given $\Phi\in C_{aa}(X,T),$
by subadditivity,  $\displaystyle
\Phi_*(\mu):=\lim_{n\to\infty} \int_X \frac{\phi_n}{n}\,
\text{d}\mu$ exists for every $\mu\in \M(X,T)$. Notice that if $\Phi$ is positive(negative),
then $\Phi_\ast(\mu)> 0(<0).$

Given a subshift of finite type $(\Sigma_A, \sigma)$ and  $\Phi\in C_{aa}(\Sigma_A,\sigma)$. If there exists  constant $D(\Phi)\ge0$ such that for any $n\ge 0,$
$$
\sup\{|\phi_n(x)-\phi_n(y)|:x|_n=y|_n\}\le D(\Phi)
$$
Then we say $\Phi$ has {\it  bounded variation} property.


\subsubsection{Thermodynamical formalism } \

Given $\Phi\in C_{aa}(\Sigma_A,\sigma)$, the topological pressure is defined as 
$$
P(\Phi):=\lim_{n\to\infty} \frac{1}{n}\ln \sum_{|w|=n}\exp(\sup_{x\in [w]}\phi_n(x)).
$$
The following extension of the classical variational principle holds:

\begin{thm}{\cite{B,M,CFH,BD}}
 Let $(\Sigma_A,\sigma)$ be a topologically mixing subshift of finite type. For any $\Phi\in C_{aa}(\Sigma_A,\sigma)$ we have 
\begin{equation}\label{supre}
P(\Phi)=\sup\{h_\mu(T)+\Phi_\ast(\mu):\mu\in \M(\Sigma_A,\sigma)\}.
\end{equation}
\end{thm}
Combining with  the monotonicity of pressure, this variational principle has the following consequence:

\begin{cor}\label{convex}
For any $\Phi, \Psi\in C_{aa}(\Sigma_A,\sigma)$, the function $Q(s):=P(\Phi+s\Psi)$ is convex on $\R.$ Consequently $Q$ is continuous. If moreover $\Psi\in C_{aa}^-(\Sigma_A,\sigma)$, then $Q(s)$ is strictly decreasing and 
$$
\lim_{s\to-\infty} Q(s)=\infty \ \ \ \text{ and }\ \ \ \lim_{s\to\infty} Q(s)=-\infty.
$$
Thus $Q(s)=0$ has a unique solution.
\end{cor}

\begin{thm}\cite{B,M} \label{gibbs}
Assume  $\Phi\in C_{aa}(\Sigma_A,\sigma)$ has bounded variation property.  Then

(i)\  There exists a unique invariant measure $\mu_\Phi$  such that 
$$
C^{-1}\le \frac{\mu_\Phi([w])}{\exp(-nP(\Phi)+\phi_n(x))}\le C \ \  (\forall x\in [w]). 
$$
$\mu_\Phi$ is called the Gibbs measure  related to $\Phi.$

(ii)\ If $\tilde \Phi\in C_{aa}(\Sigma_A,\sigma)$ and $D>0$ is a constant  such that  $\|\phi_n-\tilde\phi_n\|\le D$ for any $n\ge 0$, then $\tilde \Phi$ also has bounded variation property. Moreover $P(\Phi)=P(\tilde \Phi)$ and  $\mu_\Phi=\mu_{\tilde \Phi}$.

(iii)\ $\mu_\Phi$ is strongly mixing and  is the unique invariant measure which attains the supremum of \eqref{supre}. 

\end{thm}

If $\mu\in \M(\Sigma_A,\sigma)$ such that
$P(\Phi)=h_\mu(\sigma)+\Phi_\ast(\mu),$ then $\mu$ is called an
{\it equilibrium state} of $\Phi.$ The above proposition shows that every
$\Phi\in C_{aa}(\Sigma_A,\sigma)$ with  bounded variation property  has a unique equilibrium state.


\subsubsection{Multifractal analysis }\label{intro-multi}\

We recall some results proved in \cite{BQ}(see also \cite{BD}), which we will need in this paper.

Given a $\Psi\in C_{aa}^-(\Sigma_A,\sigma)$,
define on $(\Sigma_A,\sigma)$ a weak Gibbs metric $d_\Psi$ as 
\begin{equation}\label{weak-gibbs}
d_\Psi(x,y)=\sup_{z\in[x\wedge y]} \exp(\psi_{|x\wedge y|}(z)).
\end{equation}
This kind of metric is considered in  \cite{GP97,KS04,BQ}. In the following we will work on the metric space $(\Sigma_A,d_\Psi).$

Let $\Phi\in C_{aa}(\Sigma_A,\sigma)$ and $\Theta\in C_{aa}^-(\Sigma_A,\sigma)$. For
any $\beta\in \R$, define the level set
$$
\Lambda_{\Phi/\Theta}(\beta):=\Big\{x\in \Sigma_A:
\lim_{n\to\infty}\frac{\phi_n(x)}{\theta_n(x)} =\beta\Big\}.
$$
Since $\Theta$ is  negative,
$\Theta_\ast(\mu)<0$ for any $\mu\in \M(\Sigma_A,\sigma).$ Define
$$
L_{\Phi/\Theta}:=\Big\{\frac{\Phi_\ast(\mu)}{\Theta_\ast(\mu)} : \mu\in
\M(\Sigma_A,\sigma)\Big\}.
$$
Then $L_{\Phi/\Theta}$ is an interval.
For any $q,\beta\in \R$ we define ${\mathcal
L}_{\Phi/\Theta}(q,\beta)$ to be the unique solution
$t=t(q,\beta)$ of the equation $P(
q(\Phi-\beta\Theta)+t\Psi)=0.$  For any $\beta\in L_{\Phi/\Theta}$, define
$$
{\mathcal L}_{\Phi/\Theta}^\star(\beta)= \inf_{q\in\R}{\mathcal
L}_{\Phi/\Theta}(q,\beta).
$$

\begin{thm}\label{multifractal}\cite{BQ}
 $
 \Lambda_{\Phi/\Theta}(\beta)\ne \emptyset
\Leftrightarrow\beta\in L_{\Phi/\Theta}.
$
If $\beta\in L_{\Phi/\Theta}$, then 
$$
\dim_H \Lambda_{\Phi/\Theta}(\beta)={\mathcal L}_{\Phi/\Theta}^\star(\beta).
$$
\end{thm}
The following dimension formula for Gibbs measure is also useful.
\begin{thm}\label{exact-dim}\cite{BD}
Given $\Phi\in C_{aa}(\Sigma_A, \sigma)$ with bounded variation property. Let $\mu_\Phi$ be the related Gibbs measure. Then $\mu_\Phi$ is exact dimensional and 
\begin{equation} \label{dim-formula}
\dim_H \mu_\Phi= \frac{h_{\mu_\Phi}(\sigma)}{-\Psi_\ast(\mu)}.
\end{equation}

\end{thm}

Finally we say some words about notations. Given two positive sequences $\{a_n\}$ and $\{b_n\}$,  $a_n\lesssim b_n$ means that there exists some constant $C>0$ such that $  a_n\le C b_n$ for all $n\in \N.$ $a_n\sim b_n$ means that $a_n\lesssim b_n$ and $b_n\lesssim a_n$. $a_n \lesssim_{\gamma_1,\cdots, \gamma_k} b_n$ means that $a_n \lesssim b_n$ with the constant $C$ only depending on $\gamma_1, \cdots, \gamma_k.$


\section{A geometric lemma}\label{geo-lemma}

In this section we will prove a geometric lemma, which claims  that the ratios  of the lengths  of a gap and  the minimal band which contains it are bounded from below.  This lemma is fundamental in our  study on the metric property of the spectrum because through  it we can establish a Lipschitz equivalence between the symbolic space and the spectrum. Thus all the dimension problems of the spectrum  are completely converted to those  of  the symbolic space, where we have dynamical  tools to use.

Fix $V>20$ and $\alpha=[a_0;a_1,a_2,\cdots]$ with $a_i\le M.$
 Write 
$
{\rm Co}(\Sigma_{\alpha,V})\setminus \Sigma_{\alpha,V}=\bigcup_i G_i,
$ 
where ${\rm Co}(\Sigma_{\alpha,V})$ is the convex hull of $\Sigma_{\alpha,V}.$
Each $G_i$ is called a {\it  gap} of the spectrum. A gap $G$ is called of {\it order} $n$, if $G$ is covered by some band in $\B_n $ but not covered by any band in $\B_{n+1}.$ Denote by $\G_n$  the set of gaps of order $n$. For any $G\in \G_n$, let $B_G$ be the unique band in $\B_n$ which contains $G.$
 
 \begin{lem} \label{gap-esti}
 There exists a constant $C=C(M,V)\in (0,1)$ such that  for any gap $G\in \bigcup_{n\ge 0}\G_n$ we have 
$
|G|\ge C|B_G|.
$
\end{lem}

\proof\ Write $q:=a_{n+1}$, then $q\le M$ by the assumption. Given $B_w\in \B_n$, we study the gaps of order $n$ contained in $B_w$. If $w$ has type $t_w=I$, then there exists a unique band $B_{we_{q,q+2}}\in \B_{n+1}$ which is contained in $B_w$, thus there is no gap of order $n$ in $B_w.$

Now assume $t_w=II$. Then by \cite{R}, there exist $2q+1$ bands of order $n+1$, which are disjoint and ordered as follows:
$$
B_{we_{q,1}}<B_{we_{q,q+3}}<B_{we_{q,2}}<B_{we_{q,q+4}}<B_{we_{q,3}}<\cdots<B_{we_{q,2q+2}}<B_{we_{q,q+1}}.
$$
Thus there are $2q$ gaps of order $n$  in $B_w$. We list them as $\{G_1,\cdots,G_{2q}\}$. 

Let $S_p(x)$ be the Chebyshev polynomial defined by 
$$
S_0(x)\equiv 0,\ \  S_1(x)\equiv 1, \ \ S_{p+1}(x)=xS_p(x)-S_{p-1}(x) \ \ (p\ge 1).
$$
It is well known that 
\begin{equation}\label{cheyb}
S_p(2\cos \theta)=\frac{\sin p\theta}{\sin \theta}, \ \ \  \theta\in [0,\pi].
\end{equation}

Following \cite{FLW}, for each $p\in\N$ and $1\le l\le p$, we define 
$$
I_{p,l}:= \{2\cos \frac{l+c}{p+1}\pi: |c|\le \frac{1}{10}\ \text{ and } \  |S_{p+1}(2\cos \frac{l+c}{p+1}\pi)|\le \frac{1}{4}\}.
$$
It is seen that $I_{p,l}, l=1,\cdots,p$ are disjoint.

\noindent{\bf Claim:}  $I_{p,l}$ and $I_{p-1,s}$ are disjoint. 

\noindent$\lhd$ Fix any $x\in I_{p-1,s}$, write $x=2\cos\theta.$ Then 
\begin{equation}\label{two-re}
\frac{s-1/10}{p}\pi\le\theta  \le \frac{s+1/10}{p}\pi \ \ \text{ and }\ \ \left|\frac{\sin p\theta}{\sin\theta}\right|=|S_p(2\cos\theta)|\le \frac{1}{4}.
\end{equation}
We need to show that  $x\not \in I_{p,l}$. If otherwise, by the definition of $I_{p,l}  $ and  \eqref{two-re} we should  have 
$$
|\cos p\theta|-\left|\frac{\sin p\theta}{\sin\theta}\right|\le\left|\frac{\sin (p+1)\theta}{\sin\theta}\right|=|S_{p+1}(2\cos\theta)|\le \frac{1}{4}.
$$
Thus we have $|\cos p\theta|\le 1/2.$ On the other hand, still by \eqref{two-re} we have 
$$
(s-1/10)\pi\le p\theta\le (s+1/10)\pi.
$$
Thus 
$
|\cos p\theta|\ge \cos \frac{\pi}{10}>\frac{1}{2},
$
which is a contradiction. 
\hfill $\rhd$

By \eqref{cheyb} and the claim above,
it is easy to check that  the following intervals are disjoint subintervals of $[-2,2]$ and  ordered as 
$$
I_{p+1,1}< I_{p,1}<I_{p+1,2}<\cdots <I_{p+1,p}<I_{p,p}< I_{p+1,p+1}.
$$
   There are $2p$ gaps $\{\tilde G_1^{( p)},\cdots, \tilde G_{2p}^{(p)}\}$. Define 
$$
g( p):=\min\{|\tilde G_j^{(p)}|: j=1,\cdots,2p\}\ \ \ \text{ and }\ \ \ g:=\min\{g(1),\cdots, g(M)\}.
$$ Then $g>0$ is a constant only depending on $M.$ 

Assume $x_\ast\in B_w$ such that 
\begin{equation}\label{eq-1}
|h_w^\prime(x_\ast)||B_w|=4.
\end{equation}
By \cite{FLW} Proposition 3.1,   we have 
$$
h_w(B_{we_{q,l}})\subset I_{q+1, l}\ \ (1\le l\le q+1) \ \ \text{ and }\ \  h_w(B_{we_{q,l+q+2}})\subset I_{q, l}\ \ (1\le l\le q).
$$
Since $h_w: B_w\to [-2,2]$ is a bijection, $h_w(G_j)\supset \tilde G_j^{( q)}$ for each $j$. Fix a gap $G_j$ in $B_w$, then there exists $x_j\in G_j$ such that 
\begin{equation}\label{eq-2}
|h_w^\prime(x_j)||G_j|=|h_w(G_j)|\ge |\tilde G_j^{( q)}|.
\end{equation}
Recall that $q\le M.$ By Theorem \ref{bvar}, \eqref{eq-1}, \eqref{eq-2}, there exist a constant $c(M,V)>0$ such that 
$$
\frac{|G_j|}{|B_w|}\ge\frac{g}{4}\cdot \frac{|h_w^\prime(x_\ast)|}{|h_w^\prime(x_j)|} \ge cg=: C_1(M,V)>0.
$$ 

If $t_w=III$, the same proof as above shows that there exists a constant $C_2(M,V)>0$ such that for any gap $G$ in $B_w$ we have 
$$
\frac{|G|}{|B_w|}\ge  C_2(M,V)>0.
$$
Let $C=\min\{C_1,C_2\}$, the result follows.
\hfill $\Box$


\section{Hausdorff dimension of the dynamical subsets}\label{HDS}

From Section \ref{HDS} to Section \ref{MAOE}, we always fix $\kappa\in\N$ and   $\alpha\in \mathcal{F}_\kappa$ with continued fraction expansion $[a_0;a_1,\cdots,a_{\hat n},\kappa,\kappa,\cdots]$. In this section, we will  study the Hausdorff dimension and Hausdorff measure of $\Sigma_{\alpha,V}$ by the aid of the dynamical system $(\Omega^{(\kappa)},\sigma)$.  

We have coded $\Sigma_{\alpha,V}$ by a symbolic space $\Omega^{(\alpha)}.$ We also noticed that in general  it is hard to define a dynamic on $\Omega^{(\alpha)}.$ However since now $\alpha$ is of eventually constant type,  if  we shift $\Omega^{(\alpha)}$ $\hat n$ times, we get a subshift of finite type $\Omega^{(\kappa)}.$ It is this fact which makes it possible to use the machinery of thermodynamical formalism to study the Hausdorff dimensions of the spectrum, with the aid of   $(\Omega^{(\kappa)},\sigma)$.

At the technique level, our strategy is as follows: At first we  introduce a family of  subsets of $\Sigma_{\alpha,V},$ which we called dynamical subsets,  such that  each subset in this family can be coded by $(\Omega^{(\kappa)},\sigma)$ and the union of them are the whole spectrum.   With the aid of   $(\Omega^{(\kappa)},\sigma)$, we can obtain all the fractal  properties of  the dynamical subset. Finally we will  show that all the properties keep unchanged when we exhaust all the possible choices of subsets. Thus we get a global result for  the whole spectrum. We will finish  the  final step in Section \ref{global}.

\subsection{Dynamical subsets }\

Fix   $N=N_\alpha\ge 4+\hat n$   such that $A_\kappa^{N-\hat n}>0.$ Then the set of the last letters of words in $ \Omega^{(\alpha)}_{N}$ is $\A_\kappa.$ In other words,
$$
\{w_{N}: w\in  \Omega^{(\alpha)}_{N}\}=\A_\kappa.
$$
Since $\kappa$ is fixed, we write $e_i:=e_{\kappa,i}\in \A_\kappa$ for simplicity.
Define a subsets of $(\Omega^{(\alpha)}_{N})^{2\kappa+2}$ as 
\begin{equation}\label{def-D}
\D=\D(\alpha):= \{(w^{e_1},\cdots,w^{e_{2\kappa+2}})\in (\Omega^{(\alpha)}_{N})^{2\kappa+2}: w^{e_i}_N=e_i, i=1,\cdots, 2\kappa+2\}.
\end{equation}
We denote by $\vec{w}:=(w^{e_1},\cdots,w^{e_{2\kappa+2}})$ an element in $\D.$  Recall that  $X_w=\pi([w])$ for 
each $w\in  \Omega^{(\alpha)}_\ast.$ Given $\vec{w}\in \D$, define 
$$
\Sigma_{\vec{w}}:=\bigcup_{i=1}^{2\kappa+2} X_{w^{e_i}}.
$$
It is seen that $\Sigma_{\vec{w}}$ is made of  $2\kappa+2$  $N$-level basic sets of  $\Sigma_{\alpha,V}$ and 
$$\Sigma_{\alpha,V}=\bigcup_{\vec{w}\in \D}\Sigma_{\vec{w}}.
$$
Now we define a projection $\pi_{\vec{w}}: \Omega^{(\kappa)}\to \Sigma_{\vec{w}}$ as 
\begin{equation}\label{pi-w}
\pi_{\vec{w}}(\omega)= \pi(w^{\omega_0}\sigma\omega)=\pi(w^{\omega_0}\star\omega)
\end{equation}
(see the definition of $v\star w$ in Section \ref{coding}). It is ready to show  that $\pi_{\vec{w}}$ is a one-to-one map. We  call $\Sigma_{\vec{w}}$ a {\it dynamical subset} of $\Sigma_{\alpha,V}.$ 

We will study the dimension properties of  $\Sigma_{\vec{w}}$ at first, then by comparing  two different dynamical subsets $\Sigma_{\vec{v}}$ and $\Sigma_{\vec{w}}$, we obtain  a global picture(we will finish this in Section \ref{global}).

From now on until the end of this section  we will fix some $\vec{w}\in \D$.

\subsection{Almost additive potentials related to Lyapunov exponents }\

We will define some $\Psi\in C_{aa}^{-}(\Omega^{(\kappa)},\sigma)$, which
captures the exponential rate of the length of the  generating bands and  can be viewed as Lyapunov exponent function.
We will see that  $\Psi$ is intimately related to   the Hausdorff dimension of $\Sigma_{\vec{w}}$. 

For each $n\in \N$ and each $\omega\in\Omega^{(\kappa)}$, define 
 \begin{equation}\label{def-Psi}
 \psi_n(\omega):=\ln |B_{w^{\omega_0}\omega[1,\cdots,n]}|=\ln |B_{w^{\omega_0}\star\omega|_n}|,
 \end{equation}
where $|B_w|$ denote the length of $B_w$.

\begin{lem}  \label{bdvar}
$\Psi=\{\psi_n: n\ge 0\}\in C_{aa}^{-}(\Omega^{(\kappa)},\sigma)$ and $\Psi$ has bounded variation property.
\end{lem}

\proof\ Given $\omega\in \Omega^{(\kappa)}$ we have 
$
\psi_{n}(\omega)=\ln |B_{w^{\omega_0}\omega[1,\cdots,n]}|,$ $ \psi_{n+k}(\omega)=\ln |B_{w^{\omega_0}\omega[1,\cdots,n+k]}|$ and $\psi_k(\sigma^n\omega)=\ln |B_{w^{\omega_n}\omega[n+1,n+k]}|.
$
 By Theorem \ref{bco},  
$$
\frac{|B_{w^{\omega_0}\omega[1,\cdots,n+k]}|}{|B_{w^{\omega_0}\omega[1,\cdots,n]}|}\sim_{V,\alpha} \frac{ |B_{w^{\omega_n}\omega[n+1,n+k]}|}{ |B_{w^{\omega_n}}|}.
$$
Notice that  there are only finitely many different bands $B_w$ with $|w|=N+1$ and $N$ only depends on $\alpha$, thus we conclude  from the  above equation that 
$$
|\psi_{n+k}(\omega)-\psi_n(\omega)-\psi_k(\sigma^n\omega)|\le C(V,\alpha)=:C(\Psi).
$$
Thus $\Psi$ is almost additive.

Recall that  $N\ge 4+\hat n .$ By \eqref{band-length} we have 
\begin{equation}\label{esti-psi-n}
 (n+N)\ln c_3\le  \psi_n(\omega)\le   -n\ln 2.
\end{equation}
Thus $\Psi\in C_{aa}^-(\Omega^{(\kappa)},\sigma)$.

By the definition, if $u\in \Omega_n^{(\kappa)}$ and $\omega, \tilde\omega \in [u]$, then 
$$
\psi_n (\omega)=\ln |B_{w^{u_0}u[1,\cdots,n]}|=\psi_n(\tilde \omega).
$$
So $\Psi$ has bounded variation property with constant $D(\Psi)=0$.
\hfill$\Box$


 \subsection{Weak-Gibbs metric on $ \Omega^{(\kappa)}$}\

 Since  $\Psi\in C_{aa}^-(\Omega^{(\kappa)},\sigma),$ we can define the weak-Gibbs metric $d_\Psi$ on $\Omega^{(\kappa)}$ according to \eqref{weak-gibbs} as follows. Given $\omega,\tilde \omega\in \Omega^{(\kappa)}. $ If $\omega_0\ne \tilde \omega_0$, define $d_\Psi(\omega,\tilde \omega):={\rm diam}(\Sigma_{\alpha,V})$. If $\omega_0= \tilde \omega_0$ define 
$$
 d_\Psi(\omega,\tilde \omega):=
\sup_{\omega^\prime\in [\omega\wedge\tilde \omega]}\exp(\psi_{|\omega\wedge\tilde \omega|}(\omega^\prime))=|B_{w^{\omega_0}\star(\omega\wedge \tilde \omega)}|.
$$

 Denote by  $|\cdot|$  the standard metric on $\R$.  Then we have 

\begin{prop}\label{bi-lip}
 $\pi_{\vec{w}}: (\Omega^{(\kappa)}, d_\Psi)\to (\Sigma_{\vec{w}},|\cdot|)$ is a  bi-Lipschitz homeomorphism.
\end{prop}

\proof\ Given $\omega,\tilde \omega\in  \Omega^{(\kappa)}.$ Assume $\omega|_n=\tilde \omega|_n$ and $\omega_{n+1}\ne \omega_{n+1}.$ Then we have 
$
 d_\Psi(\omega,\tilde\omega) =|B_{w^{\omega_0}\star\omega|_n}|.
$
Write $x:=\pi_{\vec{w}}(\omega)$ and $y:=\pi_{\vec{w}}(\tilde \omega)$. It is seen that $x,y\in B_{w^{\omega_0}\star\omega|_n}$, consequently 
$$
|x-y|\le |B_{w^{\omega_0}\star\omega|_n}|=d_\Psi(\omega,\tilde\omega).
$$

On the other hand, since $\omega_{n+1}\ne \tilde \omega_{n+1},$ there is a gap $G$ of order $n+N$ which is contained in $(x,y).$ By Lemma \ref{gap-esti}, there exists a constant $C=C(\alpha,V)$ such that 
$$
|x-y|\ge |G|\ge C |B_G|=C|B_{w^{\omega_0}\star\omega|_n}|=C d_\Psi(\omega,\tilde \omega).
$$
Thus $\pi_{\vec{w}}$ is a bi-Lipschitz homeomorphism.
\hfill $\Box$

\begin{rem}
{\rm
This proposition is crucial for studying the dimensional properties of the spectrum and the density of states measure. 
Because by  this proposition, the metric property  on $(\Omega^{(\kappa)},d_\Psi)$ is the same with that on $(\Sigma_{\vec{w}},|\cdot|)$ (see for example \cite{F} chapter 2).    Thus we can convert the dimension problem of the spectrum completely to that of the symbolic space. We will use this proposition repeatedly in what follows. 
}
\end{rem}

\subsection{Bowen's formula, Hausdorff dimension and Hausdorff measure of $\Sigma_{\vec{w}}$}\

Since $\Psi \in  C_{aa}^-(\Omega^{(\kappa)},\sigma),$ by Corollary \ref{convex}, $P(s\Psi)=0$ has a unique solution $\tilde s_V=\tilde s_{V,\vec{w}}$. By  Lemma \ref{bdvar}, $\Psi$ has bounded variation, so does $\tilde s_V\Psi$. Let $m$ be the unique Gibbs measure related to $\tilde s_V\Psi$.  Then 

\begin{thm}\label{bowen-formula}
 $m\asymp \HH^{\tilde s_V}|_{\Omega^{(\kappa)}}$. Thus  $0<\HH^{\tilde s_V}( \Omega^{(\kappa)})<\infty$. 
 Consequently  $\dim_H m=\dim_H  \Omega^{(\kappa)} =\tilde s_V$ and  $\HH^{\tilde s_V}|_{\Omega^{(\kappa)}}$ is exact dimensional. Moreover 
 \begin{equation}\label{dim-haus}
 \tilde s_V= \frac{h_m(\sigma)}{-\Psi_\ast(m)}.
 \end{equation}
\end{thm}

\proof\
By Theorem \ref{gibbs} (i), there exists a constant $C>1$ such that 
for any $u\in  \Omega_n^{(\kappa)}$
$$
C^{-1}\le \frac{m( [u])}{\exp(\tilde s_V\psi_n(\omega))}\le C \ \ \  (\forall \omega\in [u] ).
$$
Together with  the definition of weak-Gibbs metric, we conclude that 
\begin{equation}\label{local}
C^{-1} {\rm diam}([u])^{\tilde s_V}\le m([u])\le C \ {\rm diam}([u])^{\tilde s_V}.
\end{equation}
Then by the definition of Hausdorff measure   it is ready to show  that $\HH^{\tilde s_V}|_{\Omega^{(\kappa)}}\asymp m$. Consequently  
 $
 0<\HH^{\tilde s_V}( \Omega^{(\kappa)})<\infty
 $
 and $\dim_H  \Omega^{(\kappa)}=\tilde s_V.$

On the other hand \eqref{local} obviously implies that $d_{m}(x)=\tilde s_V$ for all $\omega\in \Omega$. Thus $m$ and $\HH^{\tilde s_V}|_{\Omega^{(\kappa)}}$ are  exact dimensional. Finally by \eqref{dim-formula} we get 
$$
 \tilde s_V= \dim_H m=\frac{h_m(\sigma)}{-\Psi_\ast(m)}.
 $$
\hfill $\Box$

\begin{rem}\label{maximal-dim}
{\rm
The equality  $\tilde s_V=\dim_H\Omega^{(\kappa)}$ is known as Bowen's formula and \eqref{dim-haus} is kind of Young's formula. Both are very famous in the classical theories.  Since $\dim_H m$ reach the dimension of the whole space $\Omega^{(\kappa)}$, we say that $m$ is the measure of maximal dimension.  
}
\end{rem}

Now we introduce the notion of ``Gibbs type measure".

\begin{defi}\label{gibbs-type}
{\rm
Assume $X$ is a compact metric space and can be written as 
$$
X=\bigcup_{i=1}^m \bigcup_{j=1}^{n_i} X_{ij}
$$
where $X_{ij}$ are compact and disjoint. 
$\mu$ is a measure supported on $X$.  We say that $\mu$ is a {\it Gibbs type measure}, if there exists a topologically mixing  subshift of finite type $(\Sigma_A,\sigma)$ with alphabet $\A=\{1,\cdots,m\}$ such that for any choice $\tau=(\tau_1,\cdots,\tau_m), 1\le \tau_i\le n_i; 1\le i\le m$, there exist a weak Gibbs metric $d_\tau$ on $\Sigma_A, $ a Gibbs measure $\nu_\tau$ on $\Sigma_A$ and a bi-Lipschitz map $\phi^\tau:\Sigma_A\to X_\tau=\bigcup_{i=1}^m X_{i\tau_i} $  such that 
$$
\phi^\tau([i])=X_{i\tau_i}\ \ \ \text{ and }\ \ \ \phi^\tau_\ast(\nu_\tau)\asymp\mu|_{X_\tau}.
$$  
}
\end{defi}

\begin{lem}\label{exact-gibbs}
If $X$ is defined as above and $\mu$ is a Gibbs type measure supported on $X$, then $\mu$ is exact dimensional.
\end{lem}

\proof\ By Theorem \ref{exact-dim}, for each $\tau$, $\nu_\tau$ is exact dimensional. Since $\phi^\tau$ is bi-Lipschitz, we conclude that  $\mu|_{X_\tau}$ is exact dimensional. Especially $\mu|_{X_{i\tau_i}}$ is exact dimensional and the dimensions of $\mu|_{X_{i\tau_i}}$ are the same for $1\le i\le m.$ Since we can choose $\tau$ freely, we conclude that $\mu=\mu|_X$ is exact dimensional.
\hfill $\Box$

 By using Propsition \ref{bi-lip} we have the following consequence. 

\begin{thm}\label{main-Hausdorff}
$\HH^{\tilde s_V}|_{\Sigma_{\vec{w}}}$ is a Gibbs type measure, consequently $\HH^{\tilde s_V}|_{\Sigma_{\vec{w}}}$ is exact dimensional and  $0<\HH^{\tilde s_V}(\Sigma_{\vec{w}})<\infty.$ Thus $\dim_H \Sigma_{\vec{w}}=\tilde s_V$ and  
 \begin{equation}\label{dim-maxi-haus}
 \tilde s_V= \frac{h_m(\sigma)}{-\Psi_\ast(m)}.
 \end{equation}
\end{thm}

\proof\ By Proposition \ref{bi-lip}, $\pi_{\vec{w}}$ is a bi-Lipschitz homeomorphism. By general principle on Hausdorff measure(see for example \cite{F}) We conclude that 
$
(\pi_{\vec{w}})_\ast(\HH^{\tilde s_V}|_{\Omega^{(\kappa)}})\asymp \HH^{\tilde s_V}|_{\Sigma_{\vec{w}}}.
$ By Theorem \ref{bowen-formula}, we get $(\pi_{\vec{w}})_\ast(m)\asymp \HH^{\tilde s_V}|_{\Sigma_{\vec{w}}}$, thus $\HH^{\tilde s_V}|_{\Sigma_{\vec{w}}}$ is a Gibbs type measure. By Lemma \ref{exact-gibbs}, it is exact dimensional.
The other results follow from Theorem \ref{bowen-formula} and the fact that $\pi_{\vec{w}}$ is bi-Lipschitz.
\hfill $\Box$


\section{The density of states measure}\label{DSMO}

In this section we study $\NN_{\alpha,V}$. We will show that   $\NN_{\alpha,V}$ is a Markov measure and in some sense the  measure of maximal entropy. Recall that $\alpha=[a_0;a_1,\cdots,a_{\hat n};\kappa,\kappa,\cdots].$

\subsection{$\NN_{\alpha,V}$ is a  Markov measure}\

Let $H_n$ be the restriction of  $H_{\alpha,V,0}$ to the box $[1,q_{n}]$ with periodic boundary condition. Let  
$
\mathcal X_n=\{x_{n,1},\cdots, x_{n,q_n}\}
$ be the eigenvalues of $H_n$. Recall that $\sigma_n$ is defined in \eqref{sigmak}.

\begin{lem}\label{eigenvalue}
Each band in $\sigma_n$ contains exactly one value in $\mathcal X_n.$
\end{lem}

\proof \  This comes from the Bloch theory. Write $u^n=(v_1\cdots v_{q_n})^\Z$ and define $H^{(n)}=H_{u^n},$ then $\sigma_n=\sigma(H_{u^n})$ is made of  $q_n$ disjoint bands. There is another way to represent  the spectrum by using the Bloch wave.
Consider the solution of $H^{(n)}\psi_\theta=x \psi_\theta$ of the Bloch type, i.e. $\psi_\theta(m)=e^{im\theta}u(m)$ with $u(m)=u(m+q_n)$ and some $\theta\in [0,2\pi].$ When $\theta\in [0,2\pi]$ is fixed, $H^{(n)}\psi_\theta=x \psi_\theta$ has exact $q_n$ solutions. We denote the set of eigenvalues as $E_\theta.$ Then each band in $\sigma_n$ contains exact one point of $E_\theta,$ moreover  
$$
\sigma_n=\bigcup_{\theta\in [0,2\pi] } E_\theta.
$$
Now it is direct to check that $\mathcal X_n=E_0$ is the set of endpoint of each band which is  related to the phase $\theta=0.$ \hfill $\Box$

Define 
$$
\nu_n=\frac{1}{q_n}\sum_{i=1}^{q_n}\delta_{x_{n,i}}.
$$ 
It is well known that $\nu_n\to \NN_{\alpha,V}$ weakly(see for example \cite{CL}).

\begin{lem}\label{IDS}  There exist constants $C_\alpha>0$ and  
\begin{equation}\label{C-I-II-III}
C_{I}=\frac{\alpha_\kappa}{1+\alpha_\kappa^2},\ \ \ C_{II}=\frac{\alpha_\kappa^2}{1+\alpha_\kappa^2}\ \ \  \text{ and }\ \ \ C_{III}=\frac{\alpha_\kappa(\alpha_\kappa-1)}{1+\alpha_\kappa^2}
\end{equation}
such that  for any $n> \hat n$ and $B_w\in \B_n$ we have 
$$
\NN_{\alpha,V}(B_w)=C_\alpha C_{t_w} \alpha_\kappa^{-n}.
$$
\end{lem}

\proof\  By Lemma \ref{primitive}, $\hat A_\kappa$ has eigenvalues $\alpha_\kappa, -1$ and $-1/\alpha_\kappa$, then there exists invertible matrix $P$ such that $\hat A_\kappa= P\cdot {\rm diag}(\alpha_\kappa,-1,-1/\alpha_\kappa)\cdot P^{-1}.$ By \eqref{recur-q}, for $l\ge \hat n$ we have 
\begin{equation*}
\begin{pmatrix}
q_l\\
q_{l-1}
\end{pmatrix}=
\begin{pmatrix}
\kappa&1\\
1&0
\end{pmatrix}^{l-\hat n}
\begin{pmatrix}
a_{\hat n}&1\\
1&0
\end{pmatrix}\cdots 
\begin{pmatrix}
a_1&1\\
1&0
\end{pmatrix}
\begin{pmatrix}
1\\
0
\end{pmatrix}.
\end{equation*}
From this it is easy to show that there exist two constants $c_\alpha>0, d_\alpha$ such that 
\begin{equation}\label{q-l}
q_l=c_\alpha \alpha_\kappa^l+d_\alpha (-\alpha_\kappa)^{-l}.
\end{equation}

Notice that $a_k=\kappa$ for any $k\ge n$ since $n > \hat n.$
Thus for any $l>n$,  by Lemma \ref{sigma-k} and Lemma \ref{eigenvalue}, we have 
\begin{eqnarray*}
\nu_l(B_w)&=&\sum_{|u|=l+1, w\prec u}\nu_l (B_u)\\
&=&\frac{\# \{u: |u|=l+1, w\prec u, t_u=II, III\} }{q_l}\\
&=&\frac{v^{t_w}\cdot \hat A_\kappa^{l-n}\cdot v_\ast}{c_\alpha \alpha_\kappa^l+d_\alpha (-\alpha_\kappa)^{-l}}\ \ \ \  (\text{by } \eqref{q-l} \text{ and } \eqref{num})\\
&=&\frac{v^{t_w}\cdot P\cdot{\rm diag}(\alpha_\kappa^{l-n}, (-1)^{l-n}, (-\alpha_\kappa)^{n-l})\cdot P^{-1}\cdot v_\ast}{c_\alpha \alpha_\kappa^l+d_\alpha (-\alpha_\kappa)^{-l}} \\
&=&\frac{C_{t_w}\alpha_\kappa^{l-n}+C_{t_w,2} (-1)^{l-n}+C_{t_w,3} (-\alpha_\kappa)^{n-l} }{c_\alpha \alpha_\kappa^l+d_\alpha (-\alpha_\kappa)^{-l}}.
\end{eqnarray*}
Since $B_w\cap \Sigma_{\alpha,V}$ is open and closed  in $\Sigma_{\alpha,V}$ and $\nu_l\to\NN_{\alpha,V}$ weakly, 
by taking a  limit we get 
$$
\NN_{\alpha,V}(B_w)=\frac{C_{t_w}}{c_\alpha}\alpha_\kappa^{-n}=:C_\alpha C_{t_w}\alpha_\kappa^{-n}.
$$
By a simple computation, we get \eqref{C-I-II-III}.
\hfill $\Box$

Recall that we  denote $e_j:=e_{\kappa,j}$ for simplicity. Define a matrix $Q=(q_{e_ie_j})$ of order $2\kappa+2$ as 
 \begin{equation}\label{def-Q}
 q_{e_ie_j}=
 \begin{cases}
\frac{C_{t_{e_j}}}{C_{t_{e_i}}}\cdot\alpha_\kappa^{-1} & e_i\to e_j\\
 0& \text{otherwise}
 \end{cases}
 \end{equation}
 
\begin{prop}
 $Q$ is a primitive stochastic matrix. 
 Let ${p}=(p_{e_1},\cdots,p_{e_{2\kappa+2}})$ be the stationary distribution of $Q$, i.e. the unique probability vector ${p}$ such that ${p}Q={p}$, then 
\begin{equation}\label{mu-Q-e}
 p_{e_{\kappa+2}}= \frac{\alpha_\kappa}{\kappa\alpha_\kappa+2}.
\end{equation}
\end{prop}

\proof\   By the definition, $q_{e_ie_j}>0$ if and only if $A_\kappa(i,j)>0$. Since $A_\kappa$ is primitive, so is $Q$. 
Denote the $i$-th row of $Q$ by $q^i$. By the definition of $q_{e_ie_j}$ and \eqref{C-I-II-III}, we have 
\begin{equation}\label{stru-Q}
q^i=
\begin{cases}
(\underbrace{0,\cdots, 0}_{\kappa+1},1,\underbrace{0,\cdots,0}_{\kappa})& i=1,\cdots,\kappa+1\\
(\underbrace{\frac{1}{\alpha_\kappa^2},\cdots, \frac{1}{\alpha_\kappa^2}}_{\kappa+1},0,\underbrace{\frac{\alpha_\kappa-1}{\alpha_\kappa^2},\cdots,\frac{\alpha_\kappa-1}{\alpha_\kappa^2}}_{\kappa})& i=\kappa+2\\
(\underbrace{\frac{1}{\alpha_\kappa(\alpha_\kappa-1)},\cdots, \frac{1}{\alpha_\kappa(\alpha_\kappa-1)}}_{\kappa},0,0,\underbrace{\frac{1}{\alpha_\kappa},\cdots,\frac{1}{\alpha_\kappa}}_{\kappa-1},0)& i=\kappa+3,\cdots,2\kappa+2.
\end{cases}
\end{equation}
It is seen that $Q$ is a stochastic matrix.
 
We write $p_i=p_{e_i}$ temporarily.  Write $\delta:=\alpha_\kappa^{-2}$ and $\beta:=(\alpha_\kappa(\alpha_\kappa-1))^{-1}$.     Since $p=pQ$, we have 
$$
\begin{cases}
\delta p_{\kappa+2}+\beta (p_{\kappa+3}+\cdots+p_{2\kappa+2})&=p_i\ \ \ \ \ (i=1\cdots, \kappa)\\
\delta p_{\kappa+2}&=p_{\kappa+1}\\
p_1+\cdots+p_{\kappa+1}&=p_{\kappa+2}.
\end{cases}
$$
Then we    have 
$$
p_{\kappa+2}=p_1+\cdots+p_{\kappa+1}=(\kappa+1)\delta p_{\kappa+2}+\kappa\beta(p_{\kappa+3}+\cdots+p_{2\kappa+2}).
$$
Notice that $p$ is a probability vector, thus we have 
$$
1=(p_1+\cdots+p_{\kappa+1})+p_{\kappa+2}+(p_{\kappa+3}+\cdots+p_{2\kappa+2})=(2+\frac{(\alpha_\kappa-1)^2}{\alpha_\kappa})p_{\kappa+2}.
$$
Then we get $p_{\kappa+2}=\alpha_{\kappa}/(\kappa\alpha_\kappa+2).$
\hfill $\Box$

Now we have the following structure for $\NN_{\alpha,V}:$

\begin{prop}\label{markov-ids}
For any $n>\hat n$ and any $u\in \Omega^{(\alpha)}_{ n}$,  the maesure  $\NN_{\alpha,V}|_{B_u}$ is a Markov measure with transition matrix $Q$. 
 \end{prop}
 
\proof\  
For any $uw_1\cdots w_k\in \Omega_{ n+k}^{(\alpha)}$, by  using Lemma \ref{IDS} repeatedly, we have 
\begin{eqnarray*}
\NN_{\alpha,V}(X_{uw_1\cdots w_k})&=&\NN_{\alpha,V}(X_{u})\cdot\frac{\NN_{\alpha,V}(X_{uw_1})}{\NN_{\alpha,V}(X_{u})}\cdot \cdots \cdot\frac{\NN_{\alpha,V}(X_{uw_1\cdots w_k})}{\NN_{\alpha,V}(X_{uw_1\cdots w_{k-1}})}\\
&=&\NN_{\alpha,V}(X_{u})q_{w_0w_1}q_{w_1w_2}\cdots q_{w_{k-1}w_k},
\end{eqnarray*}
where $w_0$ is the last letter of $u.$ Then the result follows.
\hfill $\Box$

\subsection{$\NN_{\alpha,V}$ and the measure of maximal entropy }\

To study the dimension property of $\NN_{\alpha,V}$, we need to go back to the TDS $(\Omega^{(\kappa)},\sigma)$ again. We will establish the relation between  $\NN_{\alpha,V}$ and the measure of maximal entropy and obtain the dimension formula of $\NN_{\alpha,V}$.

Define stochastic matrix $Q$ according to \eqref{def-Q}. Let $\mu_Q$ be the unique invariant  Markov measure on the subshift $(\Omega^{(\kappa)},\sigma)$ with transition matrix $Q.$  It is well known that $\mu_Q$ is a Gibbs measure on $\Omega^{(\kappa)}$ with additive  potential   $\phi: \Omega^{(\kappa)}\to\R$ defined by 
$$
\phi(\omega):=\ln q_{\omega_0\omega_1}.
$$
Note that  $\phi\le 0$.
Write  $\Phi=(\phi_n)_{n=0}^\infty$ with $\phi_n=S_n\phi$, then $\Phi$ is a family of additive potentials.  By the definition of Gibbs measure it is not hard to compute that 
$
P(\Phi)=0.
$

\begin{thm} \label{max-entropy}
 $\mu_Q$ is the  measure of maximal entropy and  the Parry measure of the subshift $(\Omega^{(\kappa)},\sigma).$  Moreover  $\mu_Q $ is exact dimensional with 
 \begin{equation}\label{mu-Q-dim-for}
 \dim_H \mu_Q=\frac{h_{\mu_Q}(\sigma)}{-\Psi_\ast(\mu_Q)}=\frac{\ln \alpha_\kappa}{-\Psi_\ast(\mu_Q)}.   
\end{equation}
\end{thm}

\proof\  At first since $(\Omega^{(\kappa)},\sigma)$ is a topologically mixing subshift and the incidence matrix $A_\kappa$ has Perron-Frobenius eigenvalue $\alpha_\kappa,$ we have $h_{top}(\sigma)=\ln \alpha_\kappa$(see for example \cite{W}).

Recall that $p$ is the stationary distribution satisfying $pQ=p.$ Then for any $u\in \Omega_n^{(\kappa)}$, we have $\mu_Q([u])=p_{u_0}q_{u_0u_1}\cdots q_{u_{n-1}u_n}.$ For each $e\in \A_\kappa,$
fix some  $w^e\in \Omega_N^{(\alpha)}$  such that $w_N^e=e$. Then by Proposition \ref{markov-ids}, 
$$
 \mu_Q([u])=p_{u_0}q_{u_0u_1}\cdots q_{u_{n-1}u_n} =\frac{p_{u_0}}{\NN_{\alpha,V}(X_{w^{u_0}})}\NN_{\alpha,V}(X_{w^{u_0}u[1,\cdots,n]}).
$$
Then  by Lemma \ref{IDS},  for any $\omega\in \Omega^{(\kappa)},$
$$
\lim_{n\to\infty}\frac{-\ln \mu_Q([\omega|_n])}{n}=\lim_{n\to\infty}\frac{-\ln \NN_{\alpha,V}(X_{w^{\omega_{0}}\omega[1,\cdots,n]})}{n}=\ln\alpha_{\kappa}=h_{top}(\sigma).
$$
By Shannon-McMillan-Breiman Theorem, we conclude that $h_{\mu_Q}(\sigma)=h_{top}(\sigma).$ Since the measure of maximal entropy is unique, which is the so called Parry measure (see for example \cite{W} chapter 8), we conclude that $\mu_Q$ is the the measure of maximal entropy and the Parry measure  of the system $(\Omega^{(\kappa)},\sigma)$.

Since $\mu_Q$ is a Gibbs measure, by Theorem \ref{exact-dim}, $\mu_Q$ is exact dimensional and the Hausdorff dimension of $\mu_Q$ is given by \eqref{mu-Q-dim-for}. 
\hfill $\Box$

\begin{rem}
{\rm
 The proof indeed gives that for any $u\in \Omega_n^{(\kappa)}$
\begin{equation}\label{asym-mu}
\mu_Q([u])\sim \alpha_\kappa^{-n}.
\end{equation}
}
\end{rem}

We fix again a $\vec{w}\in \D$ and let $\pi_{\vec{w}}$ be defined as in \eqref{pi-w}. 
\begin{thm} \label{main-DSM}
  $\NN_{\alpha,V}|_{\Sigma_{\vec{w}}}$ is a  Gibbs type measure. It is  exact dimensional  and 
  \begin{equation}\label{NN-dim-for}
 \dim_H \NN_{\alpha,V}|_{\Sigma_{\vec{w}}}= \frac{\ln \alpha_\kappa}{-\Psi_\ast(\mu_Q)}.   
\end{equation}
  \end{thm}

\proof\  
Write $\nu_Q:=(\pi_{\vec{w}})_\ast (\mu_Q)$.
By Proposition \ref{markov-ids} we have 
$$
\nu_Q(X_{w^{u_0}u[1,\cdots,n]})=\mu_Q([u])=p_{u_0}q_{u_0u_1}\cdots q_{u_{n-1}u_n} =\frac{p_{u_0}}{\NN_{\alpha,V}(X_{w^{u_0}})}\NN_{\alpha,V}(X_{w^{u_0}u[1,\cdots,n]}).
$$
Consequently $\nu_Q\asymp \NN_{\alpha,V}|_{\Sigma_{\vec{w}}}$. Since $\mu_Q$ is a Gibbs measure, we conclude that $\NN_{\alpha,V}|_{\Sigma_{\vec{w}}}$  is a Gibbs type measure.
 Since $\pi_{\vec{w}}:\Omega^{(\kappa)}\to \Sigma_{\vec{w}}$ is bi-Lipschitz and $\nu_Q\asymp \NN_{\alpha,V}|_{\Sigma_{\vec{w}}}$, we conclude that both measures  $\nu_Q$ and $\NN_{\alpha,V}|_{\Sigma_{\vec{w}}}$ are exact dimensional and  has the same Hausdorff dimension with $\mu_Q.$
\hfill $\Box$


\section{Multifractal analysis and optimal H\"older exponent  of $\NN_{\alpha,V}$}\label{MAOE}

In this section we study the optimal H\"older exponent of $\NN_{\alpha,V}$ restricted to the dynamical subset $\Sigma_{\vec{w}}$. 
We will see that the exponent can be obtained  from the  multifractal analysis of $\NN_{\alpha,V}|_{\Sigma_{\vec{w}}}$.

At first we  conduct the multifractal analysis of $\mu_Q$, then by the bi-Lipschitz homeomorphism $\pi_{\vec{w}}$, the result is  converted  to that of $\NN_{\alpha,V}|_{\Sigma_{\vec{w}}}.$ 

\subsection{Multifractal analysis of $\mu_Q$}\

We begin with two useful lemmas:

\begin{lem}\label{two-loc-dim}
For any $\omega\in \Omega^{(\kappa)}$, we have 
\begin{equation}\label{convert}
\underline{d}_{\mu_Q}(\omega)=\liminf_{n\to\infty}\frac{\phi_n(\omega)}{\psi_n(\omega)}\ \ \ \text{ and }\ \ \ \overline{d}_{\mu_Q}(\omega)=\limsup_{n\to\infty}\frac{\phi_n(\omega)}{\psi_n(\omega)}.
\end{equation}
\end{lem}

\proof\  Fix   $\omega\in \Omega^{(\kappa)}$ and $r>0$ very small.  Assume $n$ is the unique number such that 
\begin{equation}\label{esti-r}
\exp(\psi_{nn_0}(\omega))= |B_{w^{\omega_0}\omega[1,\cdots,{nn_0}]}|<r\le|B_{w^{\omega_0}\omega[1,\cdots,{(n-1)n_0}]}|=\exp(\psi_{(n-1)n_0}(\omega)),
\end{equation}
where $n_0$ is given in Lemma \ref{band-ratio}. Thus 
$
[\omega|_{nn_0}]\subset B(\omega,r)\subset [\omega|_{(n-1)n_0}].
$
Consequently
\begin{equation}\label{esti-mu}
\mu_Q([\omega|_{nn_0}]) \le \mu_Q(B(\omega,r))\le \mu_Q([\omega|_{(n-1)n_0}]).
\end{equation}
Notice that for any $u\in \Omega_n^{(\kappa)}$ and $\omega^\prime\in [u]$, we have   
$$
\mu_Q([u])=p_{u_0}q_{u_0u_1}\cdots q_{u_{n-1}u_n}=p_{u_0}\exp(S_n\phi(\omega^\prime))=p_{u_0} \exp(\phi_n(\omega^\prime)).
$$
Combine \eqref{esti-r} and \eqref{esti-mu} we get 
$$
 \frac{\phi_{(n-1)n_0}(\omega)+\ln p_{\omega_0}}{\psi_{nn_0}(\omega)}\le \frac{\ln \mu_Q(B(\omega,r))}{\ln r}\le  \frac{\phi_{nn_0}(\omega)+\ln p_{\omega_0}}{\psi_{(n-1)n_0}(\omega)}.
$$
By \eqref{band-length},  $c_3^{N+nn_0}\le r\le 2^{2-N-(n-1)n_0}$, thus $n\to\infty$ when $r\to 0.$   By \eqref{asym-mu}, $\phi_n(\omega)\le -n\ln\alpha_\kappa+c$. By \eqref{esti-psi-n}, $\psi_n(\omega)\le -n\ln2.$
Now by using the definition of almost additive potential and  taking the upper and lower limits we get the result. 
\hfill $\Box$

\begin{lem}
There exists $d_1=d_1(\kappa)< d_2=d_2(\kappa)<0$ such that for any $\mu\in \M(\Omega^{(\kappa)},\sigma)$, we have 
\begin{equation}\label{est-Phi-mu}
d_1\le \Phi_\ast(\mu)\le d_2.
\end{equation}
\end{lem}

\proof\ Since $\Phi$ is additive and  $\mu$ is invariant, we have 
$$
\Phi_\ast(\mu)=\int_{\Omega^{(\kappa)}}\phi d\mu=\int_{\Omega^{(\kappa)}} \frac{S_3\phi}{3} d\mu.
$$ 
By the definition we get 
$$
S_3\phi(\omega)=\ln q_{\omega_0\omega_1}q_{\omega_1\omega_2}q_{\omega_2\omega_3}.
$$

We discuss two cases. At first  we assume $\kappa=1.$ By \eqref{stru-Q} we know that 
$$
q_{e_ie_j}
\begin{cases}
=1& (i,j)=(1,3); (2,3);(4,1)\\
\in (0,1)& (i,j)=(3,1);(3,2);(3,4)\\
=0& \text{otherwise}
\end{cases}
$$
 It is seen that if $e_{i_0}e_{i_1}e_{i_2}$ is admissible, then there exists at least one $j\in \{0,1,2\}$ such that $i_j=3.$   Write $q_{\min}=\min\{q_{e_3e_1},q_{e_3e_2},q_{e_3e_4}\}$ and $q_{\max}=\max\{q_{e_3e_1},q_{e_3e_2},q_{e_3e_4}\}$ then  $0<q_{\min}\le q_{\max}<1.$  Thus 
$$
3\ln q_{\min}\le S_3\phi(\omega)\le \ln q_{\max}.
$$

Next we assume $\kappa\ge 2.$ By \eqref{stru-Q},  $q_{e_ie_j}=1$ if and only if $i\le \kappa+1$ and $j=\kappa+2.$ 
Write $q_{\min}:=\min \{q_{e_{i}e_j}: e_i\to e_j;q_{e_ie_j}\ne 1\}$ and $q_{\max}:=\max \{q_{e_{i} e_j}: e_i\to e_j;q_{e_ie_j}\ne 1\}.$ Then  $0<q_{\min}\le q_{\max}<1.$    From the structure of $A_\kappa$, it is ready to see that if $e_{i_0}e_{i_1}e_{i_2}$ is admissible, then   $i_1\ne\kappa+2$ or    $i_2\ne\kappa+2$.  Thus we still have 
$$
3\ln q_{\min}\le S_3\phi(\omega)\le \ln q_{\max}.
$$
Take $d_1=\ln q_{\min}$ and $d_2=\ln q_{\max}/3$ we get the result.
\hfill $\Box$

 Consider the function $Q(q,t):=P(q\Phi+t\Psi)$. By Corollary \ref{convex}, since $\Psi\in C_{aa}^-(\Omega^{(\kappa)},\sigma),$ for each $q\in \R$ fixed, there exists a unique number $\tau(q)$ such that $Q(q,\tau(q))=0.$ Define 
 $$
 {\mathcal B}:=\{\frac{\Phi_\ast(\mu)}{\Psi_\ast(\mu)}: \mu\in \M(\Omega^{(\kappa)},\sigma)\}.
 $$
  Then ${\mathcal B}=[\beta_\ast,\beta^\ast]$ is an interval.

\begin{thm}\label{mul-ana-mu-Q}
 (i)  Define $\Lambda_\beta:=\{\omega\in \Omega^{(\kappa)}: {d}_{\mu_Q}(\omega)=\beta\}$.  Then $\Lambda_\beta\ne \emptyset $ if and only if $\beta\in {\mathcal B}.$ For any $\beta\in {\mathcal B}$
$$
\dim_H \Lambda_\beta= \tau^\ast(\beta):=\inf_{q\in \R} (\tau(q)+\beta q).
$$

(ii) There exist two constants $0<C_1(\alpha,V)\le C_2(\alpha,V)$ such that $C_1\le \beta_\ast\le \beta^\ast \le C_2.$
Moreover 
\begin{equation}\label{loc-up-low}
\beta_\ast=\inf_{\omega\in\Omega} \underline{d}_{\mu_Q}(\omega)\ \ \ \text{ and }\ \ \  \beta^\ast=\sup_{\omega\in\Omega}\overline{d}_{\mu_Q}(\omega).
\end{equation}
Thus $\beta_\ast$ is the optimal H\"older exponent of $\mu_Q.$
\end{thm}

\proof\  (i)   Recall the definition in Section \ref{intro-multi}. By \eqref{convert} we have $\Lambda_\beta=\Lambda_{\Phi/\Psi}(\beta)$ and $\tau(q)={\mathcal L}_{\Phi/\Psi}(q,\beta)-q\beta$. Thus by Theorem \ref{multifractal} we have $\Lambda_\beta\ne \emptyset$ if and only if $\beta\in {\mathcal B}$.  Moreover if $\beta\in {\mathcal B}$ we have 
$$
\dim_H \Lambda_\beta={\mathcal L}^\star_{\Phi/\Psi}(\beta)=\inf_{q\in \R}(\tau(q)+\beta q). 
$$

(ii) By \eqref{esti-psi-n}, for any invariant measure $\mu$ we have 
$$
\ln c_3\le\Psi_\ast(\mu)\le -\ln 2<0.
$$
By \eqref{est-Phi-mu}, for any invariant measure $\mu$ we have 
$$
d_1\le \Phi_\ast(\mu)\le d_2<0.
$$
From this we conclude that 
$$
0<C_1:=\frac{d_2}{\ln c_3}\le \beta_\ast\le \frac{\Phi_\ast(\mu)}{\Psi_\ast(\mu)}\le \beta^\ast\le \frac{d_1}{-\ln 2}=: C_2.
$$

Now fix any $\omega\in \Omega^{(\kappa)}.$   By Lemma \ref{two-loc-dim} we have 
$$
\underline{d}_{\mu_Q}(\omega)= \liminf_{n\to\infty} \frac{\phi_n(\omega)}{ \psi_n(\omega)}=\lim_{k\to\infty} \frac{\phi_{n_k}(\omega)}{ \psi_{n_k}(\omega)}. 
$$
By choosing a further subsequence we can further assume that 
$(\sum_{j=0}^{n_k-1}\delta_{\sigma^j\omega})/n_k\to \mu$.  Then $\mu$ is invariant and by \cite{FH} Lemma A.4(ii), we have 
$$
\lim_{k\to\infty}\frac{\phi_{n_k}(\omega)}{n_k}=\Phi_\ast(\mu) \ \ \ \text{ and }\ \ \ \lim_{k\to\infty}\frac{\psi_{n_k}(\omega)}{n_k}=\Psi_\ast(\mu).
$$
Thus $\underline{d}_{\mu_Q}(\omega)\in {\mathcal B}.$ Similarly we can show that $\overline{d}_{\mu_Q}(\omega)\in {\mathcal B}.$
On the other hand, since $\Lambda_{\beta_\ast}$ and $\Lambda_{\beta^\ast}$ are all nonempty, there exist $\omega_\ast$ and $\omega^\ast$ such that 
$d_{\mu_Q}(\omega_\ast)=\beta_\ast$ and $d_{\mu_Q}(\omega^\ast)=\beta^\ast$. Thus \eqref{loc-up-low} holds.
   \hfill $\Box$
   
   In the following we will give another expression for the optimal H\"older exponent of $\mu_Q$, which is more convenient when we study the asymptotic property. 
   
   Write $\psi_{n,\min}:= \inf\{\psi_n(\omega): \omega\in \Omega^{(\kappa)}\}.$ By the almost additivity of $\Psi$ we have 
   $$
   \psi_{n+k,\min}\ge \psi_{n,\min}+\psi_{k,\min}-C(\Psi).
   $$
   Thus $\{C(\Psi)-\psi_{n,\min}: n\ge 0\}$ form a sub-additive sequence. Then it is well known that the following limit exists 
   $$
   \lim_{n\to\infty}\frac{\psi_{n,\min}}{n}=\sup_{n\ge 0} \frac{\psi_{n,\min}-C(\Psi)}{n}=: \Psi_{\min}.
   $$
   Notice that by \eqref{esti-psi-n}, we have $\Psi_{\min}\le -\ln 2<0.$
   
\begin{prop}\label{opt-ex}
We have 
$$
\gamma_{\mu_Q}=\frac{\ln \alpha_\kappa}{-\Psi_{\min}}.
$$  
\end{prop}

\proof\ At first we show $\gamma_{\mu_Q}\ge -\ln \alpha_\kappa/\Psi_{\min}.$ It is sufficient to show that  $\underline{d}_{\mu_Q}(\omega)\ge -\ln \alpha_\kappa/\Psi_{\min}$ for any $\omega\in \Omega^{(\kappa)}.$ By the definition of $\phi_n$ and $\mu_Q$ we have  
$$
\phi_n(\omega)=\ln q_{{\omega_0}{\omega_1}}\cdots q_{{\omega_{n-1}}{\omega_n}}=\ln \frac{\mu_Q([\omega|_n])}{p_{\omega_0}q_{\omega_n\omega_{n+1}}}.
$$
By \eqref{asym-mu} we have 
\begin{equation}\label{aaa}
\lim_{n\to\infty}\frac{ \phi_n(\omega)}{n}=\lim_{n\to\infty}\frac{\ln  \mu_Q([\omega|_n])}{n}=-\ln \alpha_\kappa.
\end{equation}
On the other hand we have 
\begin{equation}\label{bbb}
\liminf_{n\to\infty} \frac{\psi_n(\omega)}{n}\ge \lim_{n\to\infty} \frac{\psi_{n,\min}}{n}=\Psi_{\min}.
\end{equation}
Combining \eqref{aaa}, \eqref{bbb} and \eqref{convert} we get
$$
\underline{d}_{\mu_Q}(\omega)=\liminf_{n\to\infty}\frac{\phi_n(\omega)}{\psi_n(\omega)}\ge \frac{\ln \alpha_\kappa}{-\Psi_{\min}}.
$$

Next we show $\gamma_{\mu_Q}\le -\ln \alpha_\kappa/\Psi_{\min}.$ It is sufficient to show that for any $\epsilon>0$ small  $\underline{d}_{\mu_Q}(\omega^\epsilon)<\ln \alpha_\kappa/(-\Psi_{\min}-\epsilon)$ for some $\omega^\epsilon\in \Omega^{(\kappa)}$. By \eqref{aaa},  it is sufficient to show that    
$\liminf_{n\to\infty}\psi_n(\omega^\epsilon)/n\le \Psi_{\min}+\epsilon$ for some $\omega^\epsilon\in \Omega^{(\kappa)}$. We find such a $\omega^\epsilon$ as follows.  Recall that the incidence matrix $A_\kappa$ is primitive, thus there exists $N_\kappa$ such that $A_\kappa^{N_\kappa-2}$ are positive.
At first take $n_0$ big enough  such that 
\begin{equation}\label{choice-n0}
\frac{\|\psi_{N_\kappa}\|_\infty}{n_0}\le \frac{\epsilon}{8}, \  \frac{C(\Psi)}{n_0}\le \frac{\epsilon}{8}, \ -4N_\kappa\Psi_{\min} \le n_0\epsilon\  \text{ and }\ \psi_{n_0,\min}\le n_0(\Psi_{\min}+\epsilon/4).
\end{equation}  
 
Let $\tilde \omega\in \Omega^{(\kappa)}$ such that $\psi_{n_0}(\tilde \omega)=\psi_{n_0,\min}.$ Since $A_\kappa^{N_\kappa-2}$ is positive, we can find $w\in \Omega_{N_\kappa-2}^{(\kappa)}$ such that both $u:=\tilde \omega|_{n_0}w$  and $w\tilde \omega|_{n_0}$ are  admissible. Thus $\omega^\epsilon:=u^\infty\in \Omega^{(\kappa)}$. Notice that $|u|=|\tilde \omega|_{n_0}|+|w|=n_0+1+N_\kappa-1=n_0+N_\kappa.$ Let  $n_1=n_0+N_\kappa.$ By the definition of $\Psi$  we have $\psi_{n_0}(\tilde\omega)=\psi_{n_0}(\omega^\epsilon).$ By almost additivity and \eqref{choice-n0} we have 
\begin{eqnarray}\label{ccc}
\psi_{n_1}(\omega^\epsilon)&\le& \psi_{n_0}(\omega^\epsilon)+\psi_{N_\kappa}(\sigma^{n_0}\omega^\epsilon)+C(\Psi)
\le \psi_{n_0}(\tilde \omega)+\frac{n_0\epsilon}{8}+\frac{n_0\epsilon}{8}\\
\nonumber&\le& n_0(\Psi_{\min}+\frac{\epsilon}{4})+\frac{n_0\epsilon}{4} \le n_1(\Psi_{\min}+3\epsilon/4).
\end{eqnarray}
Notice that by the definition of $\omega^\epsilon$ we have $\sigma^{jn_1}\omega^\epsilon=\omega^\epsilon$ for any $j\ge  0.$ Again by  almost additivity we get 
$$
\psi_{kn_1}(\omega^\epsilon)\le \sum_{j=0}^{k-1}\psi_{n_1}(\sigma^{jn_1}\omega^\epsilon)+(k-1)C(\Psi)\le k\psi_{n_1}(\omega^\epsilon)+kC(\Psi).
$$
Combining with \eqref{ccc} and \eqref{choice-n0} we conclude that 
$$
\psi_{kn_1}(\omega^\epsilon)\le kn_1(\Psi_{\min}+\epsilon).
$$
Consequently we have 
$$
\liminf_{n\to\infty}\frac{\psi_n(\omega^\epsilon)}{n}\le \liminf_{k\to\infty}\frac{\psi_{kn_1}(\omega^\epsilon)}{kn_1}\le \Psi_{\min}+\epsilon.
$$
\hfill $\Box$


\subsection{ Multifractal analysis of $\NN_{\alpha,V}|_{\Sigma_{\vec{w}}}$ }\

Through the bi-Lipschitz homeomorphism $\pi_{\vec{w}}$, we have the following: 

\begin{thm}\label{mul-ana-nn}
Let ${\mathcal B}=[\beta_\ast,\beta^\ast]$ and $\tau$ be defined as above. Define $\Lambda_\beta:=\{x\in \Sigma_{\vec{w}}: {d}_{\NN_{\alpha,V}}(x)=\beta\}$, then $\Lambda_\beta\ne\emptyset$ if and only if $\beta\in {\mathcal B}.$ For any $\beta\in {\mathcal B}$
$$
\dim_H \Lambda_\beta= \tau^\ast(\beta):=\inf_{q\in \R} (\tau(q)+\beta q).
$$
Moreover $\beta_\ast$ is the optimal H\"older exponent of $\NN_{\alpha,V}|_{\Sigma_{\vec{w}}}$ and
\begin{equation}\label{formula-optimal}
\beta_\ast=\gamma_{\NN_{\alpha,V}|_{\Sigma_{\vec{w}}}}=\frac{\ln \alpha_\kappa}{-\Psi_{\min}}.
\end{equation}
\end{thm}

\proof\ Given $\omega\in\Omega^{(\kappa)},$ write $x=\pi_{\vec{w}}(\omega).$
Since  $\pi_{\vec{w}}$ is a bi-Lipschitz homeomorphism, we have $\underline{d}_{\mu_Q}(\omega)=\underline{d}_{\NN_{\alpha,V}}(x)$ and $\overline{d}_{\mu_Q}(\omega)=\overline{d}_{\NN_{\alpha,V}}(x)$.  Then Theorem \ref{mul-ana-nn} follows from Theorem \ref{mul-ana-mu-Q} and Proposition \ref{opt-ex}.
\hfill $\Box$


\section{Global picture}\label{global}

In this section we obtain the global picture and prove Theorem \ref{main} (i), (ii) and (iii) by comparing any two different dynamical subsets $\Sigma_{\vec{v}}$ and $\Sigma_{\vec{w}}$.

\subsection{Comparison of dynamical subsets }\

Recall that for any $\alpha\in \mathcal{F}_\kappa,$  we defined $\D(\alpha)$ according to  \eqref{def-D}.  For any $\vec{w}\in \D(\alpha)$ fixed, 
we defined the potential $\Psi=\Psi_{\vec{w}}$. Let $\tilde s_{V,\vec{w}}$ be the root of $P(s\Psi_{\vec{w}})=0.$ Let $d_{\vec{w}}=d_{\Psi_{\vec{w}}}$ be the weak Gibbs metric on $\Omega^{(\kappa)}.$ Let $m_{\vec{w}}$ be the Gibbs measure with potential $\tilde s_{V,\vec{w}}\Psi_{\vec{w}}$. Let $\tilde d_{V,\vec{w}}=\dim_H^{\vec{w}} \mu_Q$ be the Hausdorff dimension of $\mu_Q$ on  the metric space $(\Omega^{(\kappa)}, d_{\vec{w}}).$  Let $\tilde \gamma_{V,\vec{w}}$ be the optimal H\"older exponent of $\mu_Q$ on the metric space $(\Omega^{(\kappa)},d_{\vec{w}}).$

Now we fix  $\alpha, \tilde \alpha \in \mathcal{F}_\kappa$. Define $\D(\alpha)$ and $\D(\tilde \alpha)$ according to \eqref{def-D}. Choose  $\vec{w}\in \D(\alpha)$ and  $\vec{v}\in \D(\tilde \alpha),$ we will  compare all the quantities related to  dynamical subsets $\Sigma_{\vec{w}}$ and $\Sigma_{\vec{v}}$.

\begin{thm}\label{comparison}
 Fix $\alpha, \tilde \alpha \in \mathcal{F}_\kappa$ and   $\vec{w}\in \D(\alpha), $ $\vec{v}\in \D(\tilde \alpha).$ Then  
$d_{\vec{v}}$ is equivalent to $d_{\vec{w}}$, $m_{\vec{v}}=m_{\vec{w}}$ and 
$$
\tilde s_{V,\vec{v}}=\tilde s_{V,\vec{w}},\ \ \ \tilde d_{V,\vec{v}}=\tilde d_{V,\vec{w}} \ \ \ \text{ and }\ \ \ \tilde \gamma_{V,\vec{v}}=\tilde \gamma_{V,\vec{w}}. 
$$
\end{thm}

\proof\
At first we show that $d_{\vec{v}}$ is equivalent to $d_{\vec{w}}$. Given $\omega,\tilde \omega\in \Omega^{(\kappa)}$ and $\omega\ne \tilde \omega.$ If $\omega_0\ne \tilde \omega_0$, then 
\begin{equation}\label{zero}
d_{\vec{v}}(\omega,\tilde \omega)={\rm diam}(\Sigma_{\tilde \alpha,V}) \ \ \ \text{ and }\ \ \ d_{\vec{w}}(\omega,\tilde \omega)={\rm diam}(\Sigma_{\alpha,V}).
\end{equation}

Now assume $\omega_0=\tilde \omega_0=e.$ Then 
$$
d_{\vec{v}}(\omega,\tilde \omega)=|B_{v^{e}\star \omega\wedge\tilde \omega}| \ \ \ \text{ and }\ \ \ d_{\vec{w}}(\omega,\tilde \omega)=|B_{w^{e}\star \omega\wedge\tilde \omega}|.
$$
 Notice that $\#\D(\alpha) (\#\D(\tilde \alpha))$ is bounded by a constant only depending on $\alpha(\tilde \alpha)$. 
 It is also clear that the partial quotients of $\alpha$ and $\tilde \alpha$ are bounded by some constant only depending on $\alpha$ and $\alpha$, since $\alpha,\tilde \alpha\in \mathcal{F}_\kappa.$
 Thus by Theorem \ref{bco} we have 
 $$
 \frac{|B_{v^{e}\star \omega\wedge\tilde \omega}|}{|B_{w^{e}\star \omega\wedge\tilde \omega}|}\sim_{(\alpha,\tilde \alpha,V)} \frac{|B_{v^{e}}|}{|B_{w^{e}}|}\sim_{(\alpha,\tilde \alpha,V)} 1.
 $$
 Together with \eqref{zero} we conclude that $d_{\vec{v}}(\omega,\tilde \omega)/d_{\vec{w}}(\omega,\tilde \omega)\sim 1.$ That is, $d_{\vec{v}}$ and $d_{\vec{w}}$ are equivalent. 
 
 As a consequence the Hausdorff dimensions and the optimal H\"older exponents  of $\mu_Q$ on $(\Omega,d_{\vec{v}})$ and  $(\Omega,d_{\vec{w}})$ are equal. That is,  $\tilde d_{V,\vec{v}}=\tilde d_{V,\vec{w}}  $ and  $\tilde \gamma_{V,\vec{v}}=\tilde \gamma_{V,\vec{w}}$.

Now we show that $\tilde s_{V,\vec{v}}=\tilde s_{V,\vec{w}}$ and  $m_{\vec{v}}=m_{\vec{w}}.$
Fix $\omega\in \Omega^{(\kappa)},$   By \eqref{def-Psi} we have 
$$
 \psi_n^{(\vec{v})}(\omega) =\ln |B_{v^{\omega_0}\omega[1,\cdots,n]}| \ \ \ \text{ and }\ \ \ \psi_n^{(\vec{w})}( \omega) =\ln |B_{w^{\omega_0}\omega[1,\cdots,n]}|.
 $$
 By Theorem \ref{bco} we have 
 $$
 \frac{|B_{v^{\omega_0}\omega[1,\cdots,n]}|}{|B_{w^{\omega_0}\omega[1,\cdots,n]}|}\sim_{(\alpha,\tilde \alpha,V)} \frac{|B_{v^{\omega_0}}|}{|B_{w^{\omega_0}}|}\sim_{(\alpha,\tilde \alpha,V)} 1.
 $$
 From this we conclude that 
 $$
 |\psi_n^{(\vec{v})}(\omega)-\psi_n^{(\vec{w})}(\omega)|\lesssim_{(\alpha,\tilde \alpha,V)} 1. 
 $$
 Now by Theorem \ref{gibbs} (ii), we conclude that $P(s\Psi_{\vec{v}})=P(s\Psi_{\vec{w}})$ for any $s\in \R$, consequently they have the same zeros. That is, $\tilde s_{V,\vec{v}}=\tilde s_{V,\vec{w}}=:\tilde s_V.$  
 Since we also have 
 $$
 |\tilde s_V\psi_n^{(\vec{v})}(\omega)-\tilde s_V\psi_n^{(\vec{w})}(\omega)|\lesssim_{(\alpha,\tilde \alpha,V)} 1,
 $$
 still by Theorem \ref{gibbs} (ii), we have 
 $m_{\vec{v}}=m_{\vec{w}}.$
  \hfill $\Box$

 \subsection{Proof of Theorem \ref{main} (i), (ii) and (iii).}\ 
 
   (i)\ 
  Fix $\alpha_\ast\in \mathcal{F}_\kappa$ and  $\vec{w}_\ast\in \D(\alpha_\ast)$, define 
  $$
   s_V(\kappa):=\tilde s_{V,\vec{w}_\ast}, \ \ \  d_V(\kappa):= \tilde d_{V,\vec{w}_\ast}\ \ \ \text{ and } \ \ \  \gamma_V(\kappa):= \tilde \gamma_{V,\vec{w}_\ast}.
  $$

 Take any $\alpha\in \mathcal{F}_\kappa,$ recall that $\Sigma_{\alpha,V}=\bigcup_{\vec{w}\in\D(\alpha)} \Sigma_{\vec{w}}.$ 
  By Theorem \ref{comparison} and Theorem \ref{main-Hausdorff}, for any $\vec{v}\in \D(\alpha)$ we have 
  $$
\dim_H \Sigma_{\vec{v}}=\tilde s_{V,\vec{v}} = \tilde s_{V,\vec{w}_\ast}=  s_V(\kappa).  
$$
Consequently 
$$
s_V(\alpha)=\dim_H \Sigma_{\alpha,V}= \dim_H \bigcup_{\vec{v}\in\D(\alpha)} \Sigma_{\vec{v}}= s_V(\kappa).
$$ 

  By Theorem \ref{comparison}, Theorem \ref{max-entropy} and Theorem \ref{main-DSM}, for any $\vec{v}\in \D(\alpha)$ we have 
  $$
\dim_H \NN_{\alpha,V}|_{\Sigma_{\vec{v}}} =\tilde d_{V,\vec{v}}=\tilde d_{V,\vec{w}_\ast}= d_V(\kappa)
$$
and $\NN_{\alpha,V}|_{\Sigma_{\vec{v}}}$ is exact dimensional.
Then by \eqref{dim-meas},  
$$
d_V(\alpha)=\dim_H \NN_{\alpha,V}= d_V(\kappa).
$$

By Theorem \ref{comparison}, Proposition \ref{opt-ex} and Theorem \ref{mul-ana-nn}, for any $\vec{v}\in \D(\alpha)$ we have 
  $$
\gamma_{\NN_{\alpha,V}|_{\Sigma_{\vec{v}}}} =\tilde \gamma_{V,\vec{v}}=\tilde \gamma_{V,\vec{w}_\ast}= \gamma_V(\kappa).
$$
Then by \eqref{optimal-H},  
$$
\gamma_V(\alpha)=\gamma_{\NN_{\alpha,V}}= \gamma_V(\kappa).
$$
Then \eqref{uniform-dim} holds.

(ii)  The result follows from   Definition \ref{gibbs-type}, Theorem \ref{main-Hausdorff} and the fact that  $\Sigma_{\alpha,V}=\bigcup_{\vec{w}\in\D(\alpha)} \Sigma_{\vec{w}}$.

(iii) The result  follows from Definition \ref{gibbs-type}, Proposition \ref{markov-ids}, Theorem \ref{main-DSM} and the fact that  $\Sigma_{\alpha,V}=\bigcup_{\vec{w}\in\D(\alpha)} \Sigma_{\vec{w}}$.
  \hfill $\Box$
  
  \begin{rem}\label{rem-optimal}
  {\rm
  
  \eqref{uniform-dim}
 has the following advantage: to compute these three quantities, 
 we can choose special element in $\mathcal{F}_\kappa$ to make the computation easier. 
 Indeed in next section, we will always pick $\alpha_\kappa\in \mathcal{F}_\kappa$ to do the computation. Moreover, due to Theorem \ref{comparison}, we can fix any $w\in \D(\alpha_\kappa),$ and compute $\gamma_V(\kappa), d_V(\kappa), s_V(\kappa)$ by \eqref{formula-optimal}, \eqref{NN-dim-for} and \eqref{dim-maxi-haus},   respectively.
  }
  \end{rem}
 

\section{Asymptotic properties  and the consequences }\label{APC}

In this section we discuss the asymptotic properties of $\gamma_V(\kappa)$, $s_V(\kappa)$ and $d_V(\kappa)$ when $V\to\infty.$ In particular,   we finish the proof of Theorem \ref{main} (iv) and (v).

By Remark \ref{rem-optimal}, in this section we always fix $\alpha=\alpha_\kappa$ and  some $w\in \D(\alpha_{\kappa}).$ Recall that we simplify $\A_\kappa=\{e_{\kappa,1}, \cdots, e_{\kappa,2\kappa+2}\}$ to $\{e_1,\cdots, e_{2\kappa+2}\}$.

\subsection{Asymptotic property of $\gamma_V(\kappa)$}\

At first we note that Lemma \ref{lem-bc} implies the following useful fact: There exists a constant $c=c_\kappa> 1$ such that  for any $w\in \Omega_n^{(\alpha_\kappa)}$,
\begin{equation} \label{length-B-w}
c^{-n} V^{-(\kappa-2)|w|_{e_{\kappa+2}}-n }\le   |B_w|\le c^n V^{-(\kappa-2)|w|_{e_{\kappa+2}}-n},
\end{equation}
where $|w|_{e_{\kappa+2}}$ stands for $\#\{1\le i\le n: w_i= e_{\kappa+2}\}.$
The proof is a direct computation by noticing that $V>20.$

\begin{prop}\label{hat-varrho-kappa}
$$
\lim_{V\to\infty} \gamma_V(1)\ln V=\frac{3}{2 }\ln \alpha_1=:\hat \varrho_1\ \ \text{ and }\ \ \lim_{V\to\infty} \gamma_V(\kappa)\ln V=\frac{2}{\kappa }\ln \alpha_\kappa=:\hat \varrho_\kappa\ \ (\kappa\ge 2).
$$
\end{prop}

\proof\ Fix some $\vec{w}\in \D(\alpha_\kappa)$ and define $\Psi $ according to \eqref{def-Psi}.  By Remark \ref{rem-optimal},  
\begin{equation}\label{for-opt}
\gamma_V(\kappa)=\frac{\ln \alpha_\kappa}{-\Psi_{\min}}.
\end{equation}
Thus we only need to estimate $\Psi_{\min}. $ Recall that $\psi_n(\omega)=\ln |B_{w^{\omega_0}\omega[1,\cdots,n]}|=\ln |B_{w^{\omega_0}\star\omega|_n}|$ and $\psi_{n,\min}=\min\{\psi_n(\omega):\omega\in \Omega^{(\kappa)}\}.$ Thus $\exp(\psi_{n,\min})$ is just the minimal length of the bands $ \{B_u: |u|=n+N_\kappa+1, w^{e_j}\prec u \text{ for some } j \}.$ 

At first we assume $\kappa=1$. Assume $u=w^{e_j} v$ with $|v|=n$.   Then by \eqref{length-B-w}
\begin{equation}\label{kappa=1}
c^{-n} V^{|v|_{e_{3}}-n}\lesssim |B_u|\lesssim c^n V^{|v|_{e_{3}}-n}
\end{equation}
Notice that $e_1e_3e_4e_1$ is admissible. 
Take $\tilde u=w^{e_4}\tilde v$ such that $|\tilde v|=n$ and  $\tilde v\prec (e_1e_3e_4)^\infty$. 
Then $|\tilde v|_{e_3}\le n/3+ 1. $ Then  \eqref{kappa=1} implies that there exists some constant $C$ (  depending on $\kappa, N_\kappa, V$) such that 
 $$
   \psi_{n,\min}\le \ln |B_{\tilde u}|\le C+n\ln c-\frac{2n}{3}\ln V.
 $$
 On the other hand for any $u=w^{e_j}v$ with $|v|=n$, by the definition of the incidence matrix $A_1$, it is ready to show that $|v|_{e_3}\ge n/3-1. $
 Then  \eqref{kappa=1} implies that there exists some constant $C^\prime$ (  depending on $\kappa, N_\kappa, V$) such that 
 $$
 \psi_{n,\min}\ge \min_{u}   \ln |B_{ u}|\ge C^\prime-n\ln c-\frac{2n}{3}\ln V.
 $$
Consequently 
 $$
-\ln c-\frac{2}{3}\ln V\le  \Psi_{\min }\le \ln c-\frac{2}{3}\ln V.
 $$
 Now by \eqref{for-opt} we conclude that $\hat \varrho_1=\frac{3}{2}\ln \alpha_1.$
 
 Next  we assume $\kappa\ge 2$. Assume $u=w^{e_j} v$ with $|v|=n$. Then by \eqref{length-B-w}
\begin{equation}\label{kappa-gt-1}
c^{-n} V^{-(\kappa-2)|v|_{e_{\kappa+2}}-n}\lesssim |B_u|\lesssim c^n V^{-(\kappa-2)|v|_{e_{\kappa+2}}-n}.
\end{equation}
Notice that $\kappa-2\ge 0$ and $e_1e_{\kappa+2}e_1$ is admissible. 
Take $\tilde u=w^{e_{\kappa+2}}\tilde v$ such that $|\tilde v|=n$ and  $\tilde v\prec (e_1e_{\kappa+2})^\infty$. 
Then $|\tilde v|_{e_{\kappa+2}}\ge n/2-1. $ Then  \eqref{kappa-gt-1} implies that there exists some constant $\tilde C$ such that 
 $$
   \psi_{n,\min}\le \ln |B_{\tilde u}|\le \tilde C+n\ln c-\frac{\kappa n}{2}\ln V.
 $$
 On the other hand for any $u=w^{e_j}v$ with $|v|=n$, since $e_{\kappa+2}e_{\kappa+2}$ is not admissible, we have  $|v|_{e_{\kappa+2}}\le n/2+1. $
 Then  \eqref{kappa-gt-1} implies that there exists some constant $\tilde C^\prime$   such that 
 $$
 \psi_{n,\min}\ge \min_{u}   \ln |B_{ u}|\ge \tilde C^\prime-n\ln c-\frac{\kappa n}{2}\ln V.
 $$
Consequently 
 $$
-\ln c-\frac{\kappa}{2}\ln V\le  \Psi_{\min }\le \ln c-\frac{\kappa}{2}\ln V.
 $$
 Now by \eqref{for-opt} we conclude that $\hat \varrho_\kappa=\frac{2}{\kappa}\ln \alpha_\kappa.$
\hfill $\Box$

\begin{rem}\label{hat-rem-varrho}
 {\rm 
 When $\kappa=1, 2$,  we have 
 $$
 \hat \varrho_1=\frac{3}{2}\ln \frac{\sqrt{5}+1}{2}\ \ \ \text{ and }\ \ \ \hat \varrho_2=\ln (\sqrt{2}+1).  
 $$
  }
  \end{rem}

\subsection{Asymptotic property of $d_V(\kappa)$}\

\begin{prop} \label{varrho-kappa}
$$
\lim_{V\to\infty} d_V(\kappa)\ln V=\frac{\kappa\alpha_\kappa+2}{2\alpha_\kappa(\alpha_\kappa-1)}\ln \alpha_\kappa=:\varrho_\kappa.
$$
\end{prop}

\proof\
By \eqref{NN-dim-for} we have  
$$
d_V(\kappa)=\dim_H \mu_Q=\frac{\ln \alpha_\kappa}{-\Psi_\ast(\mu_Q)}.
$$
  Now we study $\Psi_\ast(\mu_Q)$. Since  $\mu_Q$ is ergodic, by Kingman's sub-additive ergodic theorem,   for $\mu_Q$ a.e. $\omega\in\Omega^{(\kappa)}$, we have 
$$
-\frac{\psi_n(\omega)}{n}\to -\Psi_\ast(\mu_Q).
$$
Recall that by the definition \eqref{def-Psi}, $\psi_n(\omega)=\ln |B_{w^{\omega_0}\omega[1,\cdots,n]}|. $ By \eqref{length-B-w} we have 
\begin{equation}\label{two-side-control}
 c^{-n} V^{-(\kappa-2)\left|\omega|_n\right|_{e_{\kappa+2}}-n}\lesssim |B_{w^{\omega_0}\omega[1,\cdots,n]}| \lesssim c^n V^{-(\kappa-2)\left|\omega|_n\right|_{e_{\kappa+2}}-n}.
\end{equation}
To get the value $-\Psi_\ast(\mu_Q)$, we need to know the frequency of $e_{\kappa+2}$ in a $\mu_Q$ typical point $\omega.$ Define $\varphi(\omega)=\chi_{[e_{\kappa+2}]}(\omega).$ Since $\mu_Q$ is ergodic, by \eqref{mu-Q-e} for $\mu_Q$ a.e. $\omega\in \Omega^{(\kappa)}$ we have
$$
\frac{\#\{1\le j\le n: \omega_j=e_{\kappa+2}\}}{n}=\frac{S_n\varphi(\omega)}{n}\to \int_\Omega \varphi d\mu_Q=\mu_Q([e_{\kappa+2}])=p_{e_{\kappa+2}}=\frac{\alpha_\kappa}{\kappa\alpha_\kappa+2}.
$$
 Combining with  \eqref{two-side-control} we conclude that 
$$
\frac{2\alpha_\kappa(\alpha_\kappa-1)}{\kappa\alpha_\kappa+2}\ln V-\ln c\le -\Psi_\ast(\mu_Q)\le \frac{2\alpha_\kappa(\alpha_\kappa-1)}{\kappa\alpha_\kappa+2}\ln V+\ln c.
$$
Now combining  with the dimension formula we get the result.
 \hfill $\Box$
 
 \begin{rem}\label{rem-varrho}
 {\rm 
 When $\kappa=1,2$,  we have 
 $$
 \varrho_1=\frac{5+\sqrt{5}}{4}\ln \frac{\sqrt{5}+1}{2}\ \ \ \text{ and }\ \ \  \varrho_2=\ln (\sqrt{2}+1).
  $$
  }
  \end{rem}
   
\subsection{Asymptotic property of $s_V(\kappa)$}\

  The asymptotic properties of $s_V(\kappa)$ has been studied in \cite{LPW07, FLW,LQW}. Let us recall the result.
  
  For any $0\le  x\le1$ define
$$
\begin{array}{l}
{\mathbf R}(x) := \begin{pmatrix}
0&x^{(\kappa-1)}&0\\
(\kappa+1)x&0&\kappa x\\
\kappa x&0&(\kappa-1)x
\end{pmatrix}
\end{array}
$$
Let $\psi(x)$ be the spectral radius of ${\mathbf R}(x)$. Then it is seen that $\psi(0)=0$, $\psi(1)=\alpha_\kappa$ and $\psi(x)$ is continuous and strictly increasing. 
 Assume  $x_\kappa$  is the unique number such that $\psi(x_\kappa)=1,$ then  
  
 \begin{thm}[\cite{LPW07, FLW, LQW}]\label{rho-kappa}
 $$
 \lim_{V\to\infty} s_V(\kappa)\ln V= -\ln x_\kappa=: \rho_\kappa.
 $$
 \end{thm}
 
 In the following we will make $x_\kappa$ explicit. It is seen that for any $x>0$, the matrix ${\mathbf R}(x)$ is primitive, thus the spectral radius $\psi(x)$ of ${\mathbf R}(x)$ is the largest positive eigenvalue of ${\mathbf R}(x)$. By a direct computation, we have 
 $$
 \det (\lambda I_3-{\mathbf R}(x))= \lambda^3-(\kappa-1) x \lambda^2-(\kappa+1)x^\kappa \lambda -x^{\kappa+1}.
 $$
 Thus $x_\kappa$ is the unique number in $(0,1)$ such that 
 $$
 1-(\kappa-1)x_\kappa-(\kappa+1)x_\kappa^\kappa-x^{\kappa+1}=0.
 $$
 Write $y_\kappa=1/x_\kappa$, then $y_\kappa$ is the unique number in $(1,\infty)$ such that 
 $$
 y_\kappa^{\kappa+1}-(\kappa-1)y_\kappa^{\kappa}-(\kappa+1)y_\kappa-1=0.
 $$
 We claim that $\kappa-1<y_\kappa<\kappa$ when $\kappa\ge 3. $ Indeed define 
 $$
 F(y):=y^{\kappa+1}-(\kappa-1)y^{\kappa}-(\kappa+1)y-1,
 $$
  then  $F(\kappa-1)=-\kappa^2<0$ and $F(\kappa)=\kappa^\kappa-\kappa(\kappa+1)-1>0.$ Thus $F(y)=0$ has a root in $(\kappa-1,\kappa)$. On the other hand we know that $F(y)=0$ has only one root $y_\kappa$ in $(1,\infty)$, thus $\kappa-1<y_\kappa<\kappa.$
  
   \begin{rem}\label{rem-rho}
 {\rm 
 When $\kappa=1,2$,  by direct computation we have 
 $$
 \rho_1= \rho_2=\ln (\sqrt{2}+1).
  $$
  }
  \end{rem}

\subsection{Proof of Theorem \ref{main} (iv) and (v)}\

 At first we show (iv).
 The three  asymptotic properties have been established by Proposition \ref{hat-varrho-kappa}, Proposition \ref{varrho-kappa} and Theorem \ref{rho-kappa}.
 By Remark \ref{hat-rem-varrho}, \ref{rem-varrho} and \ref{rem-rho} we  have 
  $\hat\varrho_2=\varrho_2=\rho_2=\ln (1+\sqrt{2})$.  When $\kappa\ne 2$, it is direct to verify that $\hat \varrho_\kappa< \varrho_\kappa. $

  Now we show that   $\varrho_\kappa<\rho_\kappa$ for any $\kappa\ne 2$. 
 
 At first we claim that $\kappa<\alpha_\kappa<\kappa+1$. Indeed define $G(x)=x^2-\kappa x-1$, then $G(\kappa)=-1$ and $G(\kappa+1)=\kappa>0$, thus $G(x)=0$ has a root in $(\kappa,\kappa+1)$. On the other hand $G(x)$ has a unique positive root, which is $\alpha_\kappa$, thus we conclude that $\kappa<\alpha_\kappa<\kappa+1.$
 
Write $\delta_\kappa:= \frac{\kappa\alpha_\kappa+2}{2\alpha_\kappa(\alpha_\kappa-1)}.$ We claim that $\delta_\kappa\le 2/3$ when $\kappa\ge 8.$ Indeed for $\kappa\ge 2$ we have 
$$
\delta_\kappa= \frac{\kappa\alpha_\kappa+2}{2\alpha_\kappa(\alpha_\kappa-1)}\le \frac{\kappa(\kappa+1)+2}{2\kappa(\kappa-1)}.
$$
By a simple computation we get $\delta_\kappa\le 2/3$ for $\kappa\ge 8.$ As a result for  $\kappa\ge 8$ we have 
$$
e^{\varrho_\kappa}=\alpha_\kappa^{\delta_\kappa}\le (\kappa+1)^{2/3}.
$$
 On the other hand 
 $$
 e^{\rho_\kappa}=y_\kappa>\kappa-1.
 $$
 Thus for $\kappa\ge 8$ we have 
 $$
 e^{\varrho_\kappa}\le (\kappa+1)^{2/3}<\kappa-1<e^{\rho_\kappa}.
 $$
 That is, $\varrho_\kappa< \rho_\kappa$ for $\kappa\ge 8.$
 
 By direct computation we get    
 $\varrho_\kappa< \rho_\kappa$ for $1\le \kappa<8$ and $\kappa\ne 2$.
 Thus (iv) follows.
 
 Now we show (v).  Assume $\kappa\ne 2,$ then $\hat \varrho_\kappa<\varrho_\kappa<\rho_\kappa$. By the definition we have 
 $$
 \lim_{V\to\infty}\gamma_V(\kappa)\ln V<\lim_{V\to\infty}d_V(\kappa)\ln V<\lim_{V\to\infty}s_V(\kappa)\ln V.
 $$
 Consequently there exists $V_0(\kappa)>20$ such that for any $V\ge V_0(\kappa)$, 
 $$
\gamma_V(\kappa)\ln V< d_V(\kappa)\ln V< s_V(\kappa)\ln V.
 $$
 That is, $\gamma_V(\kappa)<d_V(\kappa)<s_V(\kappa)$ for $V\ge V_0(\kappa)$.
 \hfill $\Box$
 

\section{Appendix }\label{appendix}

In this appendix, we give another proof of the fact that  $d_V(\kappa)<s_V(\kappa)$ for $V\ge V_0(\kappa)$ when $\kappa\ne 2.$
 This proof  is more elementary and has the advantage that the constant $V_0(\kappa) $ can be estimated explicitly.

 By  \eqref{dim-maxi-haus} and \eqref{mu-Q-dim-for} we have 
$$
d_V(\kappa)=\dim_H \mu_Q=\frac{h_{\mu_Q}(\sigma)}{-\Psi_\ast(\mu_Q)}\ \ \ \text{ and }\ \ \ s_V(\kappa)=\dim_H m=\frac{h_{m}(\sigma)}{-\Psi_\ast(m)}.
$$
At first  we claim that if  $\mu_Q\ne m$, then $d_V(\kappa)<s_V(\kappa).$
Indeed by Theorem  \ref{gibbs} (iii), the unique equilibrium state of $s_V(\kappa)\Psi$ is $m$.  Since $\mu_Q\ne m$ we have 
$$
h_{\mu_Q}(\sigma)+s_V(\kappa)\Psi_\ast(\mu_Q)< h_{m}(\sigma)+s_V(\kappa)\Psi_\ast(m)=P(s_V(\kappa)\Psi)=0.
$$
Since $\Psi_\ast(\mu_Q)<0$ we conclude that  
$$
d_V(\kappa)=\dim_H \mu_Q =\frac{h_{\mu_Q}(\sigma)}{-\Psi_\ast(\mu_Q)}<s_V(\kappa).
$$

Thus we only need to study when $\mu_Q\ne m.$ For this purpose we consider two words 
$u^n=\left((II,1)_\kappa(I,1)_\kappa\right)^{3n}=(e_{\kappa+2}e_1)^{3n}$ and $\tilde u^n=\left((II,1)_\kappa(III,1)_\kappa(I,1)_\kappa\right)^{2n}=(e_{\kappa+2}e_{\kappa+3}e_1)^{2n}$. 
We will estimate respectively the following 
$$
\mu_Q([e_{1}u^n]), \ \ \mu_Q([e_{1}\tilde u^n]),\ \  m([e_{1}u^n]) \ \ \text{ and }\ \ m([e_{1}\tilde u^n]). 
$$ 

At first by \eqref{asym-mu} we have 
$$
\mu_Q([e_1u^n])\sim \alpha_\kappa^{-6n-1}\ \ \  \text{ and } \ \ \  \mu_Q([e_1\tilde u^n])\sim \alpha_\kappa^{-6n-1}.
$$
Consequently 
\begin{equation}\label{ratio-mu-Q}
\frac{\mu_Q([e_1\tilde u^n])}{\mu_Q([e_1 u^n])}\sim 1.
\end{equation}

Next we estimate  $m([e_1u^n])$ and $m([e_1\tilde u^n])$. Since $m$ is the Gibbs measure with potential $s_V(\kappa)\Psi,$ we have 
$$
m([e_1u^n])\sim |B_{w^{e_1} u^n}|^{s_V(\kappa)}\ \ \ \text{ and }\ \ \ m([e_1\tilde u^n])\sim |B_{w^{e_1}\tilde  u^n}|^{s_V(\kappa)}.
$$

Now we estimate $|B_{w^{e_1} u^n}|$ and $|B_{w^{e_1}\tilde  u^n}|$. 
At first we have 
$$
3n \le |w^{e_1} u^n|_{e_{\kappa+2}}\le N+3n\ \ \  \text{ and } 2n \le |w^{e_1} \tilde u^n|_{e_{\kappa+2}}\le N+2n.
$$

By \eqref{length-B-w}   we have 
$$
\begin{cases}
c^{-6n}V^{-3\kappa n}&\lesssim |B_{w^{e_1}u^n}|\lesssim c^{6n} V^{-3\kappa n} \\
c^{-6n}V^{-2(\kappa+1)n}&\lesssim |B_{w^{e_1}\tilde u^n}|\lesssim c^{6n} V^{-2(\kappa+1)n} \\
\end{cases}
$$

As a consequence we get 
\begin{equation}\label{ratio}
C_{V,\kappa}^n:=\left( c^{-12}V^{2-\kappa} \right)^n\lesssim \frac{|B_{w^{e_1} u^n}|}{|B_{w^{e_1}\tilde u^n}|}\lesssim \left( c^{12}V^{2-\kappa} \right)^n=:D_{V,\kappa}^n.
\end{equation}
Note that $c=c_\kappa$ is a constant only depending on $\kappa.$ 
Define $V_0(\kappa):= c_\kappa^{12}$.
  By \eqref{ratio}, it is direct to check that for $\kappa=1,$ if  $V> V_0(1)$, then   $C_{V,1}>1$; for $\kappa\ge 3$, if $V>V_0(\kappa)$, then $D_{V,\kappa}<1.$ Consequently if $V>V_0(\kappa)$, then 
\begin{eqnarray*}
\frac{m([e_1 u^n])}{m([e_1 \tilde u^n])}\sim \left(\frac{|B_{w^{e_1} u^n}|}{|B_{w^{e_1}\tilde u^n}|}\right)^{s_V(\kappa)}\gtrsim \left(C_{V,1}\right)^{ns_V(\kappa)}\to \infty ,\ \  \ (n\to\infty)& \kappa =1\\
\frac{m([e_1 u^n])}{m([e_1 \tilde u^n])}\sim \left(\frac{|B_{w^{e_1} u^n}|}{|B_{w^{e_1}\tilde u^n}|}\right)^{s_V(\kappa)}\lesssim \left(D_{V,\kappa}\right)^{ns_V(\kappa)}\to 0,\ \  \ (n\to\infty)& \kappa \ge 3. 
\end{eqnarray*}

Combine with \eqref{ratio-mu-Q} we conclude that $\mu_Q\ne m.$  Then the result follows.
\hfill $\Box$

\noindent
{\bf Acknowledgements}. The author would like to thank professor Anton Gorodetski for pointing out the references \cite{DM, Gi,Mei, Mu}, which are closely related to this paper. He also thanks professor Jean Bellissard  for pointing out that the results actually hold for more general frequencies, namely, frequencies of eventually constant type and encouraging  me to write  down this new version,  which greatly improves the previous version of the paper. Finally he thanks the referees for useful suggestions.
The author was  supported by the National Natural Science Foundation of China,   No. 11201256, No. 11371055 and No. 11431007.




\begin{thebibliography}{30}

\bibitem{BQ} J. Barral and Y. H. Qu,  \emph{  On the higher-dimensional multifractal analysis},  Discrete Contin. Dyn. Syst. 32 (2012), no. 6, 1977-1995.
 
\bibitem{B} L. Barreira,
     \emph{Nonadditive thermodynamic formalism: equilibrium and
Gibbs measures},
     Discrete Contin. Dyn. Syst., {\bf 16} (2006),
279--305.

\bibitem{BD} 
     L. Barreira and  P. Doutor,
    \emph{Almost additive
multifractal analysis},
   J. Math. Pures Appl., {\bf 92} (2009),
1--17.


 \bibitem{BIST} J. Bellissard, B. Iochum, E. Scoppola and D. Testart,
\emph{ Spectral properties of one dimensional quasi-crystals}, 
Commun. Math. Phys. {\bf 125}(1989), 527-543.



 \bibitem{C}S. Cantat, \emph{  Bers and H\'enon, Painlev\'e and Schr\"odinger}, Duke Math. J. 149 (2009), 411-460.
 
 \bibitem{Ca} M. Casdagli, \emph{ Symbolic dynamics for the renormalization map of a quasiperiodic Schr\"odinger equation},  Comm. Math. Phys. 107 (1986), no. 2, 295-318.
 
 \bibitem{CFH} 
    Y. L. Cao, D.\,J. Feng, W. Huang,
    \emph{The thermodynamic formalism
for sub-additive potentials},
    Discrete Contin. Dyn. Syst., {\bf
20} (2008), 639--657.

\bibitem{CL} R. Carmona, and J. Lacroix,  \emph{ Spectral theory of random Schršdinger operators. Probability and its Applications},  Birkh\"auser Boston, Inc., Boston, MA, 1990.


\bibitem{DEGT} D.~Damanik, M.~Embree, A.~Gorodetski, and
S.~Tcheremchantsev, \emph{ the fractal dimension of the spectrum of the
Fibonacci Hamiltonian}, Commun. Math. Phys.{\bf 280:2}(2008), 499-516.

\bibitem{DG}D. Damanik and A. Gorodetski,  \emph{ Hyperbolicity of the trace map for the weakly coupled Fibonacci
Hamiltonian}, Nonlinearity 22 (2009), 123-143.

\bibitem{DG2} D. Damanik and  A. Gorodetski, \emph{ Spectral and quantum dynamical properties of the weakly coupled Fibonacci Hamiltonian}, Commun. Math. Phys. 305 (2011), 221-277.

\bibitem{DG3}D.  Damanik and A. Gorodetski,  \emph{ The density of states measure of the weakly coupled Fibonacci Hamiltonian}, Geom. Funct. Anal. 22 (2012), no. 4, 976-989. 

 \bibitem{DG4}D.  Damanik and A. Gorodetski,  \emph{H\"older continuity of the integrated density of states for the Fibonacci Hamiltonian},  Comm. Math. Phys. 323 (2013), no. 2, 497-515.

\bibitem{DGY} D. Damanik, A. Gorodetski, W. Yessen, \emph{ The Fibonacci Hamiltonian},    arXiv:1403.7823.

\bibitem{DKL} D. Damanik, R. Killip, D. Lenz,
{\em Uniform spectral properties of one-dimensional quasicrystals,
III. $\alpha$-continuity},  Commun. Math. Phys. {\bf 212},
(2000),191-204.

\bibitem{DM}
E.  De Simone, L. Marin,  \emph{ Hyperbolicity of the trace map for a strongly coupled quasiperiodic Schr\"odinger operator},  Monatsh. Math. 163 (2011), no. 2, 211-235.

\bibitem{F} K. Falconer,   \emph{ Fractal geometry. Mathematical foundations and applications},  John Wiley \& Sons, Ltd., Chichester, 1990. 

\bibitem{Fa} K. Falconer, techniques in fractal geometry,
John Wiley\& Sons, 1997.

\bibitem{FLW} S. Fan, Q.H. Liu, Z.Y. Wen,
\emph{ Gibbs like measure for spectrum of a class of quasi-crystals},
 Ergodic Theory Dynam. Systems, {\bf 31}(2011), 1669-1695.

\bibitem{FH}  D. J. Feng and W.  Huang, \emph{  Lyapunov spectrum of asymptotically sub-additive potentials},  Comm. Math. Phys. 297 (2010), no. 1, 1-43. 

\bibitem{GM} J. B.  Garnett and D. E. Marshall, \emph{ Harmonic measure. New Mathematical Monographs, 2.} Cambridge University Press, Cambridge, 2005.


 \bibitem{GP97}
D. Gatzouras and Y.  Peres, {\em Invariant measures of full
dimension for some expanding maps},   Ergod. Th. $\&$ Dynam. Sys.,
{\bf 17} (1997),  147--167.

\bibitem{Gi} A. Girand, \emph{
Dynamical Green Functions and Discrete Schr\"odinger Operators with Potentials Generated by Primitive Invertible Substitution},
 Nonlinearity,  27 (2014) 527-543. 

\bibitem{JL} S. Jitomirskaya and Y. Last, \emph{ Power-law subordinacy and singular spectra. II. Line operators},
Commun. Math. Phys. 211 (2000), 643-658.

\bibitem{KS04} M. Kesseb\"ohmer and  B. Stratman,
{\em A multifractal formalism for growth rates and applications to
geometrically finite Kleinian groups},  Ergod. Th. $\&$ Dynam. Sys.,
{\bf 24} (2004), 141-170.


\bibitem{KKT} M. Kohmoto, L. P. Kadanoff, C. Tang, {\it Localization problem in one dimension: mapping
and escape}, {\em Phys. Rev. Lett.} {\bf 50} (1983), 1870-1872.



\bibitem{LPW07}Q.H. Liu, J. Peyri\`ere and Z.Y. Wen,
\emph{ Dimension of the spectrum of one-dimensional discrete
Schr\"odinger operators with Sturmian potentials},
 Comptes Randus Mathematique, {\bf 345:12}(2007), 667--672.

\bibitem{LQW} Q. H.  Liu, Y. H.  Qu, Z. Y.  Wen,  \emph{ The fractal dimensions of the spectrum of Sturm Hamiltonian},  Adv. Math. 257 (2014), 285-336.
 
\bibitem{LW} Q.H. Liu, Z.Y. Wen,
\emph{ Hausdorff dimension of spectrum of one-dimensional
Schr\"odinger operator with Sturmian potentials},
{ Potential Analysis} {\bf 20:1}(2004), 33--59.

\bibitem{LW05} Q.H. Liu, Z.Y. Wen,
\emph{ On dimensions of multitype Moran sets},
Math. Proc. Camb. Phyl. Soc. {\bf 139:3}(2005), 541--553.

\bibitem{Ma}  N. G. Makarov, \emph{ Fine structure of harmonic measure.}   St. Petersburg Math. J. 10 (1999), no. 2, 217-268. 

 \bibitem{Mei} May Mei, 
  \emph{   Spectra of discrete Schršdinger operators with primitive invertible substitution potentials.}  J. Math. Phys. 55 (2014), no. 8, 082701,  22 pp.
     
    
\bibitem{M} A. Mummert, \emph{ The thermodynamic formalism for almost-additive sequences}, Discrete Contin. Dyn. Syst. 16 (2006), 435-454.

\bibitem{Mu} P. Munger,
  \emph{  Frequency dependence of H\"older continuity for quasiperiodic Schr\"odinger operators, } (2013), arXiv:1310.8553.
     

\bibitem{OPRSS} S. Ostlund, R. Pandit, D. Rand, H. Schellnhuber, E. Siggia, \emph{ One-dimensional Schr\"odinger equation with an almost periodic potential}, {  Phys. Rev. Lett.} 50 (1983), 1873-1877.

 \bibitem{P} M. Pollicott,  \emph{Analyticity of dimensions for hyperbolic surface diffeomorphisms},  Proc. Amer. Math. Soc. 143 (2015), no. 8, 3465-3474.
 
\bibitem{R} L. Raymond,
{\em A constructive gap labelling for the discrete schr\"odinger
operater on a quasiperiodic chain}.(Preprint,1997)


\bibitem{Su87} A. S\"ut\"o, 
\emph{ The spectrum of a quasiperiodic Schr\"odinger operator},  Comm. Math. Phys. 111 (1987), no. 3, 409-415. 

\bibitem{Su} A. S\"ut\"o,
\emph{  Singular continuous spectrum on a Cantor set of zero Lebesgue measure for the Fibonacci Hamiltonian},  J. Stat. Phys. 56 (1989), 525-531.

\bibitem{T} M. Toda,
{\em Theory of Nonlinear Lattices},
Number 20 in Solid-State Sciences,
Springer-Verlag, second enlarged edition, 1989. Chap. 4.

 
\bibitem{W}  Walters P. \emph{ An introduction to ergodic theory. }Graduate Texts in Mathematics, 79. Springer-Verlag, New York-Berlin, 1982.
 
 \end{thebibliography}
\end{document}